\documentclass[final]{siamltex}         

\usepackage{graphicx}
\usepackage[top=26mm, bottom=26mm, left=28mm, right=28mm]{geometry}
\usepackage{array}
\usepackage[dvips]{color}
\usepackage{bm,color}
\usepackage{amsmath}
\usepackage{url}
\usepackage{algorithm}
\usepackage{algpseudocode}

\numberwithin{equation}{section}
\numberwithin{figure}{section}
\newtheorem{remark}{Remark}
\numberwithin{remark}{section}

\usepackage[T1]{fontenc}
\usepackage{amsmath}
\usepackage{graphicx}
\usepackage{amssymb}
\usepackage{bm}

\newcommand{\vct}[1]{\bm{#1}}
\newcommand{\mtx}[1]{\mathsf{#1}}

\newcommand{\vtwo}[2]{\left[\begin{array}{c} #1 \\ #2 \end{array}\right]}

\newcommand{\mtwo}[4]{\left[\begin{array}{cc} #1 & #2 \\ #3 & #4  \end{array}\right]}

\newcommand{\bi}{\begin{itemize}}
\newcommand{\ei}{\end{itemize}}
\newcommand{\ben}{\begin{enumerate}}
\newcommand{\een}{\end{enumerate}}
\newcommand{\be}{\begin{equation}}
\newcommand{\ee}{\end{equation}}
\newcommand{\bea}{\begin{eqnarray}} 
\newcommand{\eea}{\end{eqnarray}}
\newcommand{\ba}{\begin{align}} 
\newcommand{\ea}{\end{align}}
\newcommand{\bse}{\begin{subequations}} 
\newcommand{\ese}{\end{subequations}}
\newcommand{\bc}{\begin{center}}
\newcommand{\ec}{\end{center}}
\newcommand{\bfi}{\begin{figure}}
\newcommand{\efi}{\end{figure}}

\newcommand{\tbox}[1]{{\mbox{\tiny #1}}}
\newcommand{\mbf}[1]{{\mathbf #1}}
\newcommand{\pO}{{\partial\Omega}}
\newcommand{\half}{\mbox{\small $\frac{1}{2}$}}
\newcommand{\uu}{\mbf{u}}                       
\newcommand{\vv}{\mbf{v}}                       
\newcommand{\xx}{\mbf{x}}                       
\newcommand{\yy}{\mbf{y}}                       
\newcommand{\rr}{\mbf{r}}                       
\newcommand{\nn}{\mbf{n}}                       
\newcommand{\dd}{\mbf{d}}                       
\newcommand{\nx}{\nn^\xx}                       
\newcommand{\ny}{\nn^\yy}                       
\newcommand{\dx}{d_\xx}                          
\newcommand{\dy}{d_\yy}                          

\newcommand{\pdrive}{p_\tbox{drive}}        
\newcommand{\nr}{^\tbox{near}}           
\newcommand{\RP}{R_P}      
\newcommand{\tU}{\bm\tau_U}    
\newcommand{\tD}{\bm\tau_D}
\newcommand{\jump}[1]{[\![#1]\!]_\gamma}
\newcommand{\um}{\mathbf{u}}

\begin{document} 

\title{A fast algorithm for simulating multiphase flows through periodic geometries of arbitrary shape
}

\author{
Gary R. Marple\thanks{Department of Mathematics, 
University of Michigan, Ann Arbor, MI, 48109, USA. {\em email:} {\tt gmarple@umich.edu}.}
\and 
Alex Barnett\thanks{Department of Mathematics, 
Dartmouth College, Hanover, NH, 03755, USA. {\em email:} {\tt ahb@math.dartmouth.edu}.}%
\and
Adrianna Gillman\thanks{Computational and Applied Mathematics, 
Rice University, Houston, TX,  77005, USA. {\em email:} {\tt adrianna.gillman@rice.edu}.}%
\and Shravan Veerapaneni\thanks{Department of Mathematics, 
University of Michigan, Ann Arbor, MI, 48109, USA. {\em email:} {\tt shravan@umich.edu}.}}
\date{\today}
\maketitle

\begin{abstract}
This paper presents a new boundary integral equation (BIE)
method for simulating particulate and multiphase flows through periodic channels of arbitrary 
smooth shape in two dimensions. The authors consider a particular system---multiple 
vesicles suspended in a periodic channel of arbitrary shape---to describe 
the numerical method and test its performance. 
Rather than relying on the periodic Green's function as classical BIE
methods do, the method
combines the free-space Green's function
with a small auxiliary basis, and 
imposes periodicity as an extra linear condition.
As a result, we can exploit
existing free-space solver libraries, quadratures, and fast algorithms,
and handle a large number of vesicles in a geometrically complex channel.
Spectral accuracy in space is achieved using the periodic trapezoid
rule and product quadratures, while a first-order semi-implicit scheme evolves particles by treating the vesicle-channel  
interactions explicitly. New constraint-correction formulas are introduced that preserve reduced areas of vesicles,
independent of the number of time steps taken.  
By using two types of fast algorithms, (i) the fast multipole method
(FMM) for the computation of the vesicle-vesicle and the vesicle-channel hydrodynamic 
interaction, and (ii) a fast direct solver for the BIE on 
the fixed channel geometry, the computational cost is reduced to ${\cal O}(N)$ per 
time step where $N$ is the spatial discretization size. Moreover, the direct solver inverts the wall BIE operator at $t = 0$, stores its compressed representation and applies it at every time step to evolve the vesicle positions, leading to dramatic cost savings compared to classical approaches.  
Numerical experiments illustrate that a simulation with $N=128{,}000$ can 
be evolved in less than a minute per time step on a laptop.
\end{abstract}


\section{Introduction}
%
%
%
Suspensions of rigid and/or deformable particles in viscous fluids flowing through confined geometries are ubiquitous in natural and engineering systems. Examples include drop, bubble, vesicle, swimmer, and red blood cell (RBC) suspensions. Understanding the spatial distribution of such particles in confined flows is crucial in a wide range of applications including targeted drug delivery \cite{Huang2010}, enhanced oil recovery \cite{olbricht1996pore}, and microfluidics for cell sorting and separation \cite{jin2014technologies}. 
In several of these applications, the long-time behavior of the suspension is sought. For example: {\em What is the optimal size and shape of targeted drug carriers that maximizes their ability to reach the vascular walls escaping from flowing RBCs} \cite{Huang2010, Decuzzi2010}? {\em What is the optimal design of a microfluidic device that differentially separates circulating tumor cells from blood cells} \cite{Russom2009}? More generally, one is interested in estimating the rheological properties of a given particulate suspension in an applied flow, electric, or magnetic fields. A common mathematical construct that is employed in such a scenario is the periodicity of flow at the inlet and the outlet. Therefore, the natural computational problem that arises is to solve for the transient dynamics of a particulate flow through a confined periodic channel driven either by pressure difference or other stimuli. 

In this work, we consider periodization algorithms for {\em vesicle suspensions} in confined flows. Vesicles---often considered as mimics for biological cells, especially RBCs \cite{vlahovska2013flow}---are comprised of bilipid membranes enclosing a viscous fluid and their diameter is typically less than 10$\mu m$. The membrane mechanics is modeled by the Helfrich energy \cite{helfrich1973elastic} combined with a local inextensibility constraint. At the length scale of the vesicles, the Reynolds number is extremely small, therefore, the Stokes equations are employed to model the fluid interior and exterior to the vesicles. The suspension dynamics of this system is governed by the nonlinear membrane forces, the vesicle-vesicle and vesicle-channel non-local hydrodynamic interactions, and the applied flow boundary conditions. 

\begin{figure}[t]
\centering
\includegraphics[width=\textwidth]{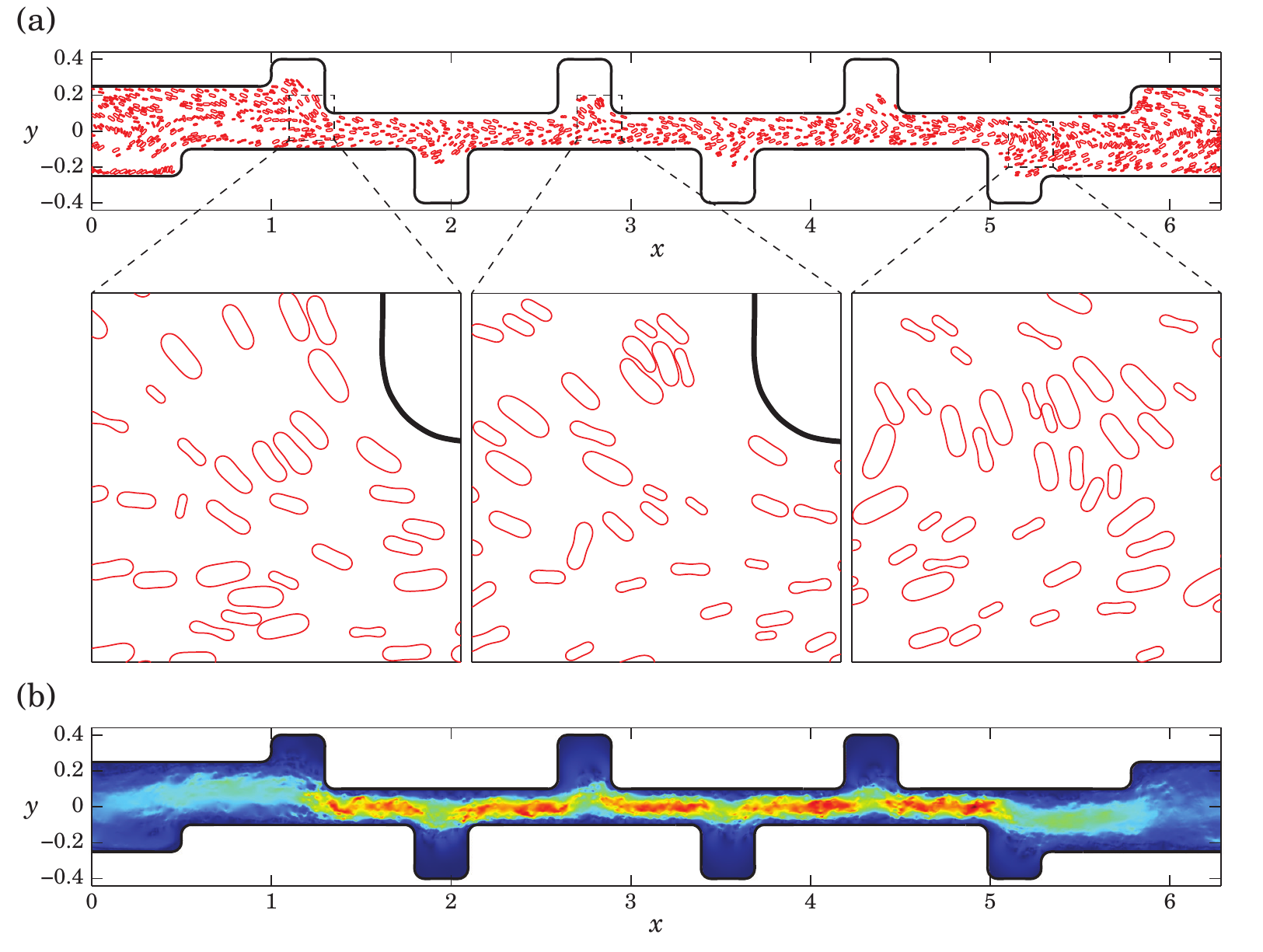}
\caption{(a) Snapshot from a simulation of 1{,}005 vesicles flowing through an arbitrary-shaped periodic channel. We used 64 discretization points per vesicle and 32,000 points each for the top and bottom walls. The vesicle-vesicle and vesicle-channel hydrodynamic interactions are computed via the Stokes FMM, and the new fast direct solver (Sec. \ref{s:fds}) is used to solve the channel BIEs. We used the close evaluation scheme of \cite{lsc2d} for the vesicle-to-vesicle and vesicle-to-channel interactions (but not for channel-to-vesicle interactions due to its inapplicability).  This simulation took 52 seconds per time step on a laptop with a 2.4 GHz dual-core Intel Core i5 processor and 8 GB of RAM. 
(b) Plot of the velocity magnitude (red indicates high and blue indicates low) corresponding to the disturbance field generated by the vesicles (obtained by subtracting the pressure-driven ``empty pipe'' flow from the total velocity field).}
\label{fig:largesim}
\end{figure}

Pioneered by Youngren and Acrivos \cite{Youngren75, Youngren76}, BIE methods are widely used for particulate and other interfacial flows \cite{pozrikidis01interfacial}. Their advantages over grid- and mesh-based methods are well-known: reduction in dimensionality, ease of achieving high-order accuracy, and availability of highly scalable fast algorithms. The classical approach for incorporating periodic boundary conditions within the BIE framework is to replace the free-space Green's function with one that satisfies the periodicity condition. This can be expressed as an infinite sum of source images. The double-layer potential defined on an open curve $\Gamma$ which defines one period of a channel wall with lattice vector ${\bf d}$, for instance, can be written as
\begin{equation}
{\bf u}({\bf x}) \,= \, {\! \sum_{{n} \in\mathbb{Z}} \,\,\int_\Gamma D({\bf x},{\bf y}+ n{\bf d}) \, \bm{\tau} ({\bf y})\, ds_{\bf y}}, \label{periodicSum} \end{equation}
where {\bf u} is the fluid velocity, $D$ is the Stokes free-space double-layer kernel, and $\bm{\tau}$ is the density function defined on $\Gamma$. Classical algorithms, such as the {\em Ewald summation} \cite{Ewald21,hasimoto1959periodic, pozrikidis1996computation} for accelerating the $N-$body calculation that arises from discretizing \eqref{periodicSum}, use a partition of unity to split the discrete sum into rapidly converging sums for the {\em nearby} and {\em distant} interactions, handling the latter in the spectral domain.
The local interactions are $\mathcal{O}(N)$ in number by construction,
leaving distant interactions which
can be evaluated accurately at the $N$
targets by combining local interpolations onto a regular grid
with the fast Fourier transform (FFT),
a technique named particle-mesh Ewald \cite{dardenPME}.
Recently, an accurate variant called spectral Ewald has been
developed for
particulate flows \cite{lindbo2010spectrally, af2014fast, wang2015spectral}
in periodic geometries.

Although such FFT-based methods are widely used, owing to their ease of implementation, they suffer from several drawbacks.
Firstly, the FFT introduces a ${\mathcal O}(N \log N)$ complexity,
and although the constants are rather small for FFT methods, the scalability of communication costs on multicore architectures is suboptimal (see \cite{gholami2014fft} for a detailed discussion). Secondly, the lack of spatial adaptivity makes them somewhat impractical for problems with multi-scale physics. 
Finally, the ``gridding'' required is expensive, becoming even more so
in three dimensions. 
One of the main goals of this paper is to introduce a simple alternative algorithm that is ${\mathcal O}(N)$, and, since it exploits existing fast algorithms,
overcomes many of these limitations which can be quite restrictive for constrained geometries that have local features (see Fig. \ref{fig:largesim}).
  
{\em Synopsis of the new method.} The proposed periodizing integral equation formulation is based on 
the ideas introduced in Barnett--Greengard \cite{qpsc} for quasi-periodic scattering problems. 
It uses direct free-space summation for the nearest-neighbor periodic images, whereas the flow 
field due to the distant images is captured using an auxiliary basis comprised of a small 
number of ``proxy'' sources. Periodic boundary conditions are imposed in an extended linear 
system (ELS) that determines both the wall layer densities (to enforce no-slip boundary 
condition on the channel) and the proxy source strengths.
Although one block of this ELS is rectangular and ill-conditioned,
its pseudo-inverse is rapid to compute, allowing
accuracy close to machine precision \cite{mlqp,acper}.
The disturbance velocities and the 
hydrodynamic stresses due to the presence of vesicles enter the right-hand side of the ELS 
so that the combined flow field is periodic from channel inlet to outlet and vanishes on the walls. 
The proposed integral formulation is versatile in handling the imposed flow boundary conditions: applied pressure-drop 
across the channel, or imposed slip on the channel (e.g., to model electroosmotic flows),
simply modify the right-hand side of the ELS.
The scheme can handle various dimensions of periodicity and easily
extends to 3D \cite{acper}.

The proposed particulate flow solver is based on the work of Veerapaneni et al. \cite{ves2d}. It 
employs a semi-implicit time-stepping scheme to overcome the numerical stiffness associated with 
the integro-differential equations governing the vesicle evolution. A spectral (Fourier) basis 
is used to represent the vesicle and channel boundaries. The required spatial derivatives are 
computed via spectral differentiation and the singular integrals are also computed with spectral 
accuracy using the product quadrature rule given in Kress \cite{kress91}. The FMM is used to accelerate 
the computation of the vesicle-vesicle hydrodynamic interactions. A simple correction term is introduced in the local inextensibility constraint applied at every time step. It eliminates error accumulation over long time-periods that usually leads to numerical instabilities.

Since the channel geometry remains fixed,
its BIE linear system may be inverted once and for all.
We use a direct solver to precompute its inverse;
the channel wall densities can then be determined at 
each time step for the small cost of a matrix-vector multiply. Since the matrix associated with wall-wall interactions
has off-diagonal low-rank structure, the
use of a hierarchical fast direct solver
reduces the cost of both the precomputation,
and of a solve with each new right-hand side to ${\mathcal O} (N)$, where $N$ is the number of points 
on the channel walls.
Crucially, the cost involved with each new solve is very small when compared
against the standard combination of an iterative solver and an FMM:
for example, Table~\ref{tab:fastcompare} shows that, in our setting,
the former is {\em three orders of magnitude} faster than the latter.
The fast direct solver for the ELS is a Stokes version of the periodic Helmholtz scattering solver
of Gillman--Barnett \cite{qps}, with the extra complication that the channel walls form a 
continuous interface (in \cite{qps} it was possible to isolate the obstacles).  The continuous 
interface demands a new compression scheme, presented
in Sec.~\ref{s:fds}.
While fast direct solvers have received much attention recently, 
the authors believe that this work is the first to apply them to particulate flows.


{\em Advantages.} The conspicuous advantage of our method is the simplicity that comes from the use of free-space kernels (as opposed to lattice sums or particle-mesh Ewald) and a pressure drop condition that is applied directly. More importantly, this feature allows us to use 
state-of-the-art high-order quadratures for the singular \cite{LIE} and nearly-singular integrals \cite{helsing_close, lsc2d} and open-source FMM implementations \cite{CMCL}. 
The periodization scheme itself is spectrally accurate in terms of the number of proxy points.  Numerical experiments illustrate that 
the same number of proxy points are required to achieve a specified accuracy independent of the complexity of the geometry.

{\em Limitations.} In this paper, we restrict our attention to two-dimensional problems. Although the periodization scheme itself is straightforward to extend to three dimensions, other components of the numerical method, such as the quadratures and the direct solver for surfaces in three dimensions, require more work. We assume that there are no holes in the given periodic domain and that the viscosities of the fluids interior and exterior to vesicle membranes are the same. Both of these assumptions are merely for simplicity of exposition (e.g., see \cite{rahimian10} to relax the latter assumption and the completed double-layer formulation \cite{Pozrikidis:BIM} for the former). 

Our recent algorithm to evaluate the nearly-singular integrals with spectral accuracy \cite{lsc2d} is applied to evaluate the vesicle-vesicle hydrodynamic interactions when they are close to each other. However, we do not apply any close evaluation schemes for the channel-to-vesicle interactions (\cite{lsc2d} cannot be applied directly to this setting). This limitation prohibits us from performing simulations of tightly-packed suspensions. One possible remedy is to switch the discretization of the channel from a global basis (Fourier) to local chunk-based schemes, for which close evaluation schemes for Stokes potentials already exist \cite{ojala2014accurate}.     

\subsection{Related Work} 
%

There is a large body of literature on BIE methods for periodic Stokesian flows of rigid and deformable bodies. A few examples include \cite{zick1982stokes, larson1986microscopic, fan1998completed, greengard2004integral} for fixed particles as in a porous medium, \cite{Loewenberg96, Zinchenko00, freund2007lmm, zhao10, af2014fast} for particulate suspensions, and particularly \cite{zhao2013dynamics, thiebaud2013rheology} for vesicle flows. However, almost all studies focus on flows through either simple geometries, such as flat channels or cylinders, or without any constraining walls (i.e., one-, two-, or three-periodic systems in free-space).
The work of Greengard--Kropinski \cite{greengard2004integral}
uses an intrinsically two-dimensional
complex-variable formulation to periodize a FMM-based solver
for a fixed doubly-periodic geometry in $\mathcal{O}(N)$ cost.
However, the imposition of pressure-drop conditions in their scheme
is a subtle matter involving non-convergent lattice sums.
In contrast, in the present work, such conditions are applied
simply and directly and the cost remains $\mathcal{O}(N)$.
The only work that we are aware of where vesicles inhabit an arbitrary periodic geometry is that by Zhao et al \cite{zhao10} where capsules flowing through deformed cylinders were simulated. The constraining geometry is embedded in a box and a Green's function that satisfies periodic conditions on the box is used. One of the drawbacks of this approach is that a pressure drop cannot be imposed directly but is determined from the mean flow. Furthermore, a large amount of auxiliary data might be introduced because of the embedding in the case of geometries that have multi-scale spatial features. 
The authors would like to point out that to date, there are no reported results on particulate flows through complex periodic geometries such as those shown in Fig. \ref{fig:largesim}.

Fast direct solvers, such as ${\mathcal H}$-matrix, HBS, HSS, and HODLR 
(\cite{2009_xia_superfast,2007_shiv_sheng,2010_xia,2010_borm_book,2004_borm_hackbusch,m2011_1D_survey,HODLR}) which utilized hierarchical 
low-rank compression of off-diagonal blocks, are naturally applicable to solving the linear systems arising from the 
discretization of boundary integral equations, thanks to the smoothly decaying property of Green's functions for 
far interactions.  Since the construction of the fast direct solver dominates the computational cost, it is not 
beneficial to build a solver for the evolving geometries.  Instead, the proposed method constructs a fast 
direct solver for the fixed constrained walls, then reuses this precomputed solver at each time step to evolve the 
vesicles.  The computational cost for both steps scales linearly with the number of discretization points. The cost for each time step is much reduced compared to an iterative FMM-based solve.
While this manuscript describes an HBS solver (see \cite{m2011_1D_survey} or Section \ref{s:fds}) for simplicity of presentation, alternative $\mathcal{O}(N)$ inversion techniques can be 
seamlessly substituted in.
%

\subsection{Outline} 
This manuscript begins by describing the periodic solver for the steady Stokes equation in a channel given velocity boundary conditions (Section 2). The numerical scheme for the evolution of vesicles and the treatment of the coupling between the vesicle and channel BIEs is presented in Section 3. The fast ELS solver for finding
the channel densities and the proxy strengths is presented in Section 4. 
The accuracy, stability, and computational complexity of the method will be illustrated via numerical experiments in Section 5. Finally, the manuscript concludes with a summary and statement of future work in Section 6.


\section{Periodization scheme}
\label{s:per}

\subsection{Preliminaries: Stokes potentials and the non-periodic BVP}

We first define the standard kernels and boundary integral operators used
\cite{Ladyzhenskaya} (for our 2D case see \cite[Sec. 2.2, 2.3]{HW}).
Let $\mu>0$ be the fluid viscosity, a scalar constant.
Let $\uu(\xx)=(u_1(\xx),u_2(\xx))$
be the velocity field and $p(\xx)$ the scalar
pressure field for $\xx=(x_1,x_2)\in\mathbb{R}^2$. The pair $(\uu,p)$
is a solution to the Stokes equations if
\bea
-\mu\Delta\uu + \nabla p &=&0
\label{force}
\\
\nabla\cdot\uu &=& 0.
\label{incomp}
\eea
These express force balance and incompressibility, respectively.

The Stokes single-layer kernel (stokeslet)
from source point $\yy$ to target point $\xx$
has tensor components
\begin{equation}
S_{ij}({\bf x},{\bf y})=\frac{1}{4\pi\mu}\left(\delta_{ij}\log{\frac{1}{r}}+\frac{r_ir_j}{r^2}\right),\hspace{12pt}i,j=1,2,
\label{S}
\end{equation}
where $\rr:=\xx-\yy$, $r:=\|\rr\|$, and $\delta_{ij}$ is the Kronecker delta.
Given a density (vector function) $\bm\tau$ on a source curve $\Gamma$,
the single-layer representation for velocity is then
$\uu = {\cal S}_\Gamma\bm\tau$, i.e.,
\be
\uu(\xx) = ({\cal S}_\Gamma \bm\tau)(\xx) :=
\int_\Gamma S(\xx,\yy) \bm\tau(\yy) ds_\yy.
\label{SLP}
\ee
The associated pressure function is
\be
p(\xx) = ({\cal Q}_\Gamma \bm\tau)(\xx) :=
\int_\Gamma Q(\xx,\yy) \bm\tau(\yy) ds_\yy
\quad \mbox{ where }
Q_j(\xx,\yy) = \frac{1}{2\pi}\frac{r_j}{r^2}
~,\qquad j=1,2.
\label{Q}
\ee
For the double-layer velocity representation
$\uu = {\cal D}_\Gamma\bm\tau$, we have, using
$\ny$ the surface normal at each point on the
source curve $\Gamma$,
\be
\uu(\xx) = ({\cal D}_\Gamma \bm\tau)(\xx) :=
\int_\Gamma D(\xx,\yy) \bm\tau(\yy) ds_\yy
\quad \mbox{ where }
D_{ij}(\xx,\yy) = \frac{1}{\pi}\frac{r_ir_j}{r^2}\frac{\rr\cdot\ny}{r^2}
~,\qquad i,j=1,2.
\label{D}
\ee
We write its associated pressure function as
\be
p(\xx) = ({\cal P}_\Gamma \bm\tau)(\xx) :=
\int_\Gamma P(\xx,\yy) \bm\tau(\yy) ds_\yy
\quad \mbox{ where }
P_j(\xx,\yy) = \frac{\mu}{\pi}\left( -\frac{\ny_j}{r^2} + 2 \rr\cdot\ny\frac{r_j}{r^4}\right)
~,\qquad j=1,2.
\label{P}
\ee

We use the notation $D_{\Gamma',\Gamma}$ to indicate
the double-layer
boundary integral operator from source curve $\Gamma$ to target
$\Gamma'$, i.e., $D_{\Gamma', \Gamma}\bm\tau = ({\cal D}_\Gamma\bm\tau)|_{\Gamma'}$.
If the target and source curves are the same ($\Gamma'=\Gamma$)
then $D_{\Gamma,\Gamma}$ is to be taken in the principal value sense
and has a smooth kernel for smooth $\Gamma$.
We have, for $\Gamma$ a $C^2$-smooth curve and any $\bm\tau\in C(\Gamma)$,
the jump relation
\be
\lim_{h\to0^+}({\cal D}_\Gamma\bm\tau)(\xx - h\nx)
= (-\half I + D_{\Gamma,\Gamma})\bm\tau
\label{JR}
\ee
for the interior limit of velocity. Here, $I$ is the $2\times2$ identity tensor.
The non-periodic prototype BVP that we will need to solve is that the
pair $(\uu,p)$ satisfies the Stokes equations
in a bounded domain $\Omega$ for given velocity (Dirichlet) data
$\uu = \vv$ on its boundary $\pO$.
To solve this problem,
we insert the double-layer representation $\uu = {\cal D}_\pO \bm\tau$
into \eqref{JR} to get the 2nd-kind boundary integral equation (BIE) on $\pO$:
\be
(-\half I + D_{\pO,\pO})\bm\tau = \vv
~.
\label{BIE}
\ee
This BVP, and the resulting BIE,
has one consistency condition, $\int_\pO \vv \cdot \nn ds_\yy = 0$,
and null-space of dimension one corresponding to adding a constant to $p$ \cite{HW}.

Since we will also need to impose {\em traction} (Neumann)
matching conditions, we need the
traction on a target curve due to the above representations.
Given a function pair $(\uu,p)$, the Cauchy stress tensor at any point
has entries
\be
\sigma_{ij}(\uu,p):= -\delta_{ij} p + \mu(\partial_i u_j + \partial_j u_i)
~,\hspace{12pt}i,j=1,2.
\label{sigma}
\ee
The hydrodynamic
traction of this pair, i.e., the force vector per unit length applied to the fluid
at a surface point with outward unit normal $\nn$,
has components
\be
T_i(\uu,p) := \sigma_{ij}(\uu,p)n_j =
-pn_i + \mu(\partial_i u_j + \partial_j u_i) n_j
~, \hspace{12pt}i=1,2,
\label{Tup}
\ee
where summation over $j$ is implied.
Applying \eqref{Tup} to the $(\uu(\xx),p(\xx))$
pair, due to the single-layer velocity
\eqref{S} and pressure \eqref{Q} kernel (with fixed source point $\yy$),
gives the single-layer traction kernel
\be
K_{ik}(\xx,\yy) = \sigma_{ij}(S_{jk}(\cdot,\yy),Q_k(\cdot,\yy)) (\xx) \nx_j
= -\frac{1}{\pi}\frac{r_ir_k}{r^2}\frac{\rr\cdot\nx}{r^2}
~,\qquad i,k=1,2,
\label{K}
\ee
which we abbreviate by $K$.
Likewise, applying \eqref{Tup} to the double-layer pair \eqref{D} and \eqref{P}
gives, after a somewhat involved calculation
(eg \cite[(5.27)]{Liu09book}),
the double-layer traction kernel tensor
\bea
T_{ik}(\xx,\yy) &=& \sigma_{ij}(D_{jk}(\cdot,\yy),P_k(\cdot,\yy)) (\xx) \nx_j
~,\qquad i,k=1,2
\nonumber \\
&=&
\frac{\mu}{\pi}\left[
\left(\frac{\ny\cdot\nx}{r^2}- 8\dx\dy\right)\frac{r_ir_k}{r^2}
+\dx\dy\delta_{ik} +\frac{\nx_i\ny_k}{r^2}
+\dx\frac{r_k \ny_i}{r^2} + \dy \frac{r_i \nx_k}{r^2}
\right],
\label{T}
\eea
where we defined the target and source ``dipole functions''
$$\dx = \dx(\xx,\yy) := (\rr\cdot\ny)/r^2
~,\hspace{1in}
\dy = \dy(\xx,\yy) := (\rr\cdot\nx)/r^2,
$$
respectively.
The use of the symbol $T$ to mean the traction operator
and the double-layer traction kernel will be clear by context.
The hypersingular
boundary integral operator for the traction of the double-layer
from source curve $\Gamma$ to target $\Gamma'$
we call $T_{\Gamma',\Gamma}$.
To clarify,
$$
(T_{\Gamma',\Gamma}\bm\tau)(\xx) = T({\cal D}_\Gamma\bm\tau, {\cal P}_\Gamma\bm\tau)(\xx) = \int_\Gamma T(\xx,\yy)\bm\tau(\yy) ds_\yy
~, \qquad \xx\in\Gamma',
$$
where in the final expression, the kernel $T$ has tensor components \eqref{T}.


\begin{figure}
\includegraphics[width=\textwidth]{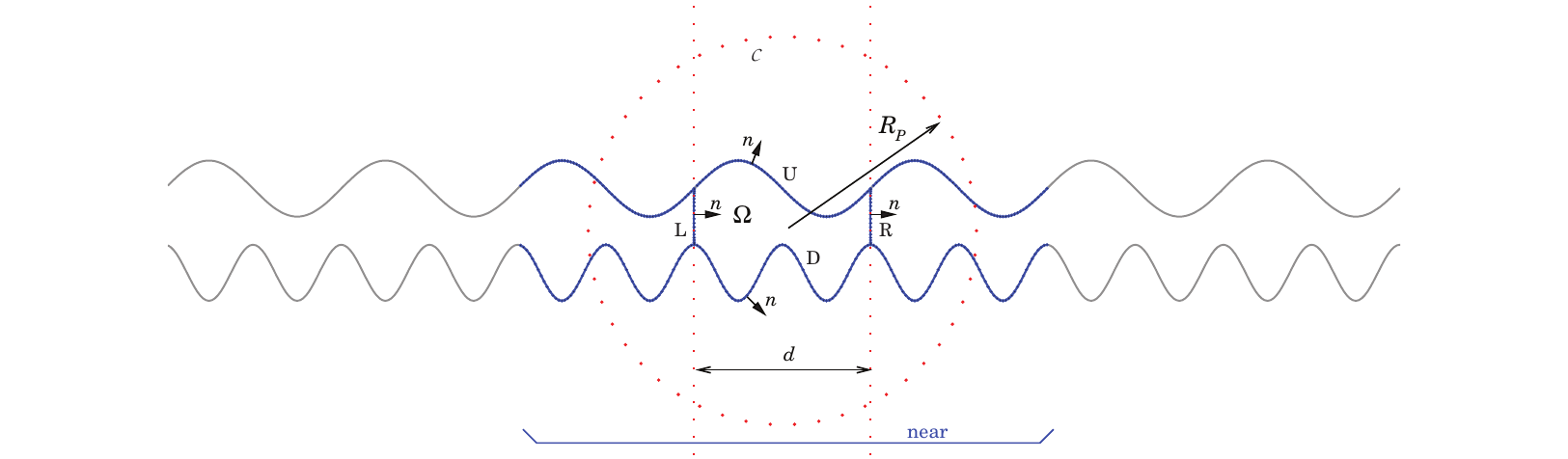}%
\caption{Geometry for periodization scheme of Section \ref{s:per}. The grey shows the infinite
periodic pipe, and the blue dots the quadrature nodes for the
central domain and its near neighbors. The central domain $\Omega$
is bounded by $\Gamma = U \cup D$, and side walls $L$ and $R$,
and has the normal senses shown.
The proxy points (red) lie on circle $\mathcal C$ of radius $R_p$.}
\label{fig:setup}
\end{figure}

\subsection{Dirichlet problem in a driven periodic pipe}

We consider a single unit cell $\Omega$
confined by one period of the walls $U$ above and $D$ below.
The full periodic pipe
domain is then
$\Omega_\Lambda := \{\xx\in\mathbb{R}^2:\xx+n\dd\in\Omega,
n\in\mathbb{Z}\}$, where $\dd = (d,0)$ is the lattice vector with period $d$.
See Fig.~\ref{fig:setup}.

It is conceptually simplest to begin with
the following strictly-periodic ``empty pipe'' problem.
Given periodic
velocity (Dirichet) data $\vv_U$ and $\vv_D$ on the up and down walls,
find a solution $(\uu,p)$ in $\Omega_\Lambda$ that
is periodic up to a constant pressure driving per period, i.e.,
\bea
(\uu,p) & &\mbox{ Stokes in } \Omega_\Lambda
\label{perBVP1}
\\
\uu &=& \vv_U \mbox{ on } U
\\
\uu &=& \vv_D \mbox{ on } D
\\
\uu(\xx+\dd) - \uu(\xx) &=& 0, \quad \xx\in\Omega_\Lambda
\\
p(\xx+\dd)- p(\xx) &=& \pdrive, \quad \xx\in\Omega_\Lambda.
\label{perBVPf}
\eea
The consistency condition on the data is
$\int_U \vv_U\cdot\nn ds_\yy + \int_D \vv_D\cdot\nn ds_\yy = 0$, and the nullity 1,
as in the non-periodic case.
In our application, the data $\vv_D$, $\vv_U$ will be
(minus) the flow velocity induced by a periodized set of
vesicles inside $\Omega_\Lambda$.

A standard approach for solving this BVP in the
strictly-periodic case $\pdrive=0$
would be to sum the double-layer kernel over all periodic copies in order to obtain the periodized version of the kernel:
\begin{equation}
D^P({\bf x},{\bf y})=\sum_{n\in\mathbb{Z}}{ D({\bf x},{\bf y}+n{\bf d})}.
\label{ksum}
\end{equation}
The representation is then, using $\Gamma = U \cup D$ to indicate the
one period of the upper and lower walls,
\begin{equation}
{\bf u}({\bf x})=(\mathcal{D}^P_\Gamma\bm{\tau})({\bf x})=\int_U{\! D^P({\bf x},{\bf y})\bm{\tau}_U({\bf y})\, ds_{\bf y}}+\int_D{\! D^P({\bf x},{\bf y})\bm{\tau}_D({\bf y})\, ds_{\bf y}}.
\end{equation}
By analogy with \eqref{BIE},
the density $\bm{\tau}=[\bm{\tau}_U;\bm{\tau}_D]$ could be
determined by solving the 2nd-kind integral equation
\begin{equation}
\left(-\half I+\mathcal{D}^P_{\Gamma,\Gamma}\right)\bm{\tau}=\vv
\end{equation}
with $\vv=[\vv_U;\vv_D]$.

\begin{remark}
If $\pdrive \neq 0$, this approach could also be used
after substracting from $\vv$ velocity data from the Poiseuille flow
$\uu(\xx) = (\frac{1}{2}\alpha x_2^2,0)$, $p(\xx) = \alpha \mu x_1$, with
$\alpha = \pdrive/(\mu d)$.
The result is a strictly-periodic BVP with $\pdrive=0$ and modified data $\vv$.
However, we will find the following approach much more convenient.
\end{remark}

Instead, we reformulate the BVP on the single unit cell $\Omega$,
introducing a left side wall $L$ and right side wall $R = L + \dd$
(Fig.~\ref{fig:setup}).
Note that, given a periodic pipe $\Omega_\Lambda$, the choice of where to place
the wall to subdivide the unit cell is arbitrary.
We choose them to be vertical for convenience.
Furthermore we relax the periodicity condition on the wall velocity data
$\vv_U$, $\vv_D$,
and impose between $L$ and $R$ periodicity conditions for velocity and traction
with given arbitrary mismatch $\mbf{g}_u$, $\mbf{g}_T$,
that we call the ``discrepancies'' \cite{qpsc}.
Thus,
\bea
(\uu,p) & &\mbox{ Stokes in } \Omega
\label{boxPDE}
\\
\uu &=& \vv_U \mbox{ on } U
\label{boxU}
\\
\uu &=& \vv_D \mbox{ on } D
\label{boxD}
\\
\uu_R - \uu_L &=& \mbf{g}_u
\label{boxgu}
\\
T(\uu,p)_R - T(\uu,p)_L &=& \mbf{g}_T.
\label{boxgT}
\eea
By unique continuation from Cauchy data, if $\vv_U$ and $\vv_D$ are periodic,
and we choose $\mbf{g}_u \equiv \mbf{0}$ and $ \mbf{g}_T = \pdrive\nn$,
where $\nn$ here indicates the normal $(1,0)$ on the $L$ and $R$ walls,
then the above BVP is equivalent to
\eqref{perBVP1}--\eqref{perBVPf}.
The special case $\vv_U\equiv\vv_D\equiv\mbf{0}$ gives pressure-driven flow in a
periodic pipe, free of vesicles.
In the general case, there is still a consistency condition on the data.
The advantage of the above (non-periodic) BVP in the single unit cell
is that the data may be induced by a sum over vesicles which includes
only the nearest images and that pressure driving is incorporated naturally.

To solve \eqref{boxPDE}--\eqref{boxgT}, we use a kernel containing
only the {\em near-field} images, plus a small auxiliary
basis for smooth Stokes solutions in $\Omega$ to account for the effect
of the infinite number of {\em far-field} images.
For the latter, we use the  ``method of fundamental solutions'' (MFS) basis
\cite{Bo85,mfs}
of stokeslets with sources
lying on a circular contour ${\cal C}$ enclosing $\Omega$.
(These are also known as ``proxy points'' \cite{m2011_1D_survey}.)
To be precise, the velocity representation is
\be
\uu \;=\; {\cal D}_\Gamma\nr\bm\tau + \sum_{m=1}^M \mbf{c}_m \phi_m,
\label{urep}
\ee
where
\begin{equation}
({\cal D}\nr_\Gamma\bm\tau)({\bf x})
:=
\sum_{|n|\leq 1}{\int_U{\! D({\bf x},{\bf y}+n{\bf d})\bm{\tau}_U({\bf y})\, ds_{\bf y}}}+\sum_{|n|\leq 1}{\int_D{\! D({\bf x},{\bf y}+n{\bf d})\bm{\tau}_D({\bf y})\, ds_{\bf y}}}
\label{Dnr}
\end{equation}
is a sum over free-space kernels living on the walls in the central unit cell and its two near neighbors.
The second term contains basis functions $\phi_m$ that satisfy the Stokes equation in the physical domain living in the central unit cell. The basis $\{\phi_m\}$ needs to accurately represent any field due to the ``far'' periodic copies (i.e. those indexed $\ldots,-3,-2,2,3,\ldots$). The source points
$\{\yy_m\}_{m=1}^M$ are
equispaced on a circle of sufficiently large radius $\RP$ centered on the central unit
cell, and
\be
\phi_m({\bf x}) = S({\bf x},{\bf y}_m), \qquad m=1,\ldots,M
\label{phi}
\ee
is the stokeslet at the $p$th proxy source.
Each coefficient $\mbf{c}_m \in \mathbb{R}^2$, totalling $2M$ unknowns.
This may be viewed as approximating a single-layer density lying on the circle which is able to represent inside any field due to sources lying outside.
Since the sources are distant from $\Omega$, the convergence is exponential
with a rapid rate; we only need $M={\cal O}(1)$ (typically less than $10^2$)
{\em independent of the complexity of the domain or the number of quadrature
points needed to accurately represent it.}

The pressure representation corresponding to \eqref{urep} is
(summing \eqref{P} in the same fashion),
\be
p = {\cal P}_\Gamma\nr\bm\tau + \sum_{m=1}^M \mbf{c}_m \varphi_m,
\qquad
\mbox{ where }
\varphi_m({\bf x}) = Q({\bf x},{\bf y}_m), \quad m=1,\ldots,M.
\label{prep}
\ee
Thus, for any density $\bm\tau$ and coefficients $\{\mbf{c}_m\}$,
$(\uu,p)$ solves the Stokes equations in $\Omega$.

\begin{remark}
There are constraints on the radius $\RP$:
larger $\RP$ allows for more rapid error convergence with respect to $M$,
but if $\RP$ is larger than $3d/2$,
then the circle encloses some image sources and
the size of the coefficients $\mbf{c}_m$ grow exponentially large,
resulting in catastrophic cancellation.
Hence, we fix $\RP = d$ in this study.
\end{remark}

Constructing a linear system is now simply a matter of inserting
the representation \eqref{urep} into each of the conditions
\eqref{boxU}--\eqref{boxgT}, which we now do.
Imposing the velocity data on $U$ and $D$
using the jump relation \eqref{JR} (which only affects the $n=0$ term
in \eqref{Dnr}) gives two coupled boundary integral-algebraic equations,
\bea
(-\half I + D\nr_{U,U})\tU + D\nr_{U,D}\tD
+ \sum_{m=1}^M \phi_m |_U \mbf{c}_m 
&=&
\vv_U \qquad \mbox{ on } U
\label{BIEU}
\\
D\nr_{D,U}\tU + (-\half I + D\nr_{D,D})\tD
+ \sum_{m=1}^M \phi_m |_D \mbf{c}_m 
&=&
\vv_D \qquad \mbox{ on } D~,
\label{BIED}
\eea
which we may summarize as
$$
A \bm\tau + B {\bm c} = {\bm v}
~.
$$
Imposing periodicity in matching velocity and traction data
\eqref{boxgu}--\eqref{boxgT} gives,
after noticing cancellations of all of the close wall-wall interactions,
\bea
(D_{R,U-\dd}-D_{L,U+\dd})\tU + (D_{R,D-\dd}-D_{L,D+\dd})\tD +
\sum_{m=1}^M (\phi_m |_R-\phi_m |_L) \mbf{c}_m
&=&
\mbf{g}_u
\label{BIEgu}
\\
(T_{R,U-\dd}-T_{L,U+\dd})\tU + (T_{R,D-\dd}-T_{L,D+\dd})\tD +
\sum_{m=1}^M \bigl(T(\phi_m,\varphi_m) |_R-T(\phi_m,\varphi_m)|_L\bigr) \mbf{c}_m
&=&
\mbf{g}_T
\label{BIEgT}
\eea
which we summarize as
$$
C \bm\tau + Q {\bm c} = {\bm g}
~.
$$
Thus, the four coupled boundary integral-algebraic equations
\eqref{BIEU}--\eqref{BIEgT}
may be stacked in pairs and written in a block form
\begin{equation}
\left[ \begin{array}{cc}
A & B\\
C & Q\\
\end{array}\right]
\left[ \begin{array}{c}
{\bm \tau}\\
{\bm c}\\
\end{array}\right]=
\left[ \begin{array}{c}
{\bm v}\\
{\bm g}\\
\end{array} \right].
\label{bsys}
\end{equation}
The roles of the block matrices are as follows:
$A$ is a 2nd-kind operator mapping wall densities to ($U$, $D$) wall velocities,
$B$ maps auxiliary coefficients to wall velocities,
$C$ maps wall densities to their discrepancies
in the periodicity conditions,
and $Q$ maps auxiliary coefficients to their discrepancies.
The functions ${\bm v}$ and ${\bm g}$ contain the Dirichlet and discrepancy data.

\subsection{Discretization of the linear system}

The coupled integral-algebraic equations \eqref{bsys} need to be discretized:
this is performed in a standard fashion using
$N$ quadrature nodes $\{\xx_{U,j}\}_{j=1}^N$ on $U$,
$N$ nodes $\{\xx_{D,j}\}_{j=1}^N$ on $D$,
and $K$ nodes $\{\xx_{L,j}\}_{j=1}^K$ on $L$ (the nodes on $R$ being
those on $L$ displaced by $\dd$).
The nodes on $U$ and $D$ are generated using the periodic trapezoid rule
applied to a smooth parametrization of the curves,
while the $L$ and $R$ wall ``collocation''
nodes are chosen to be Gauss--Legendre
in the vertical coordinate (no weights are needed for these nodes)
\cite[Ch.~9]{NA}.

The Nystr\"om method \cite[Sec.~12.2]{LIE} is used to discretize $A$.
For instance, given the quadrature weights $\{w_{U,j}\}_{j=1}^N$ on $U$,
the matrix discretization of the $1,1$ tensor block of $D_{U,U}$ has elements
\be
\mtx{D}_{ij} = \left\{ \begin{array}{ll}
D_{11}(\xx_{U,i},\xx_{U,j}) w_{U,j}, & i\neq j \\
-\frac{\kappa(\xx_{U,j})}{2\pi} (\mbf{t}_1(\xx_{U,j}))^2 w_{U,j}, & i=j
\end{array}
\right.,
\ee
which uses the diagonal limit
$\lim_{\yy\to\xx}D_{ij}(\xx,\yy) = -\frac{\kappa(\xx)}{2\pi} \mbf{t}_i(\xx) \mbf{t}_j(\xx)$,
where $\kappa$ is the curvature of the boundary,
and $\mbf{t}$ the unit tangent vector.
Other blocks are filled similarly.
The result is to replace \eqref{bsys} by a discrete linear system of
similar structure, which will be solved with a fast direct solver described 
in Section~\ref{s:fds}.
For more details on a similar periodic discretization scheme, see \cite{mlqp}.




\section{Application to particulate flows}
\label{s:appl}



Now we describe the application of the above periodic BVP solution
to vesicle flow simulations, extending
the scheme of Veerapaneni et al \cite{ves2d} to
periodic flow of vesicle suspensions in rigid pipe-like geometries.
We begin with just a single (periodized) vesicle with its boundary 
$\gamma$ lying in the periodic unit cell $\Omega$.
We solve the steady-state flow problem given forces on the vesicle
and then use this solver within the time-stepping scheme.

\subsection{Solving for quasi-static fluid flow given the interfacial forces}
\label{s:uresp}

For simplicity, we consider the case without viscosity contrast.
Let ${\bf x}(s)$ parametrize the vesicle membrane
$\gamma$ according to arc-length $s$.
The membrane generates forces on the fluid due to bending ${\bf f}_b=-\kappa_B{\bf x}_{ssss}$ and tension ${\bf f}_\sigma=\left(\sigma {\bf x}_s\right)_s$, where $\kappa_B$ is the bending modulus and $\sigma$ is the tension.
The total force at each point on $\gamma$ is then $\mbf{f} = \mbf{f}_b + \mbf{f}_\sigma$. Stress balance and no-slip condition at the membrane-fluid interface imply the jump conditions $\jump{T(\uu, p)} = \mbf{f}$ and $\jump{\uu} = 0$ respectively, where $\jump{\cdot}$ denotes the jump across $\gamma$. In the case of an isolated vesicle in free-space and assuming $\mbf{f}$ is known, the representation for the fluid velocity $\uu = {\cal S}_\gamma\mbf{f}$ and the pressure $p = {\cal Q}_\gamma\mbf{f}$ satisfies these jump conditions as well as the Stokes equations in the bulk.


Our goal in this section is, given only the forces $\mbf{f}$ on a
periodized vesicle $\gamma+n\dd$, $n\in\mathbb{Z}$,
to solve for the fluid flow $\uu$ which results
in the periodic confining geometry $\Omega_\Lambda$,
with no-slip boundary conditions and given pressure drop $\pdrive$ across the channel.
The basic idea is to write $\uu$ as a sum of the ``imposed''
flow that the vesicle would generate in an infinite fluid,
plus a ``response'' flow due to the confining geometry $\Omega$.
(This is the same concept as the incident and scattered wave in
scattering theory \cite{coltonkress}.)
A standard approach in the case $\pdrive=0$
might be to use a periodic imposed flow
$\sum_{n\in\mathbb{Z}} {\cal S}_{\gamma+n\dd}\mbf{f}$
and for the response to solve the periodic BVP
\eqref{perBVP1}--\eqref{perBVPf} with
velocity data given by the negative of the imposed flow measured
on $U$ and $D$.
The sum of imposed and response flows then would meet our goal.
However, this approach has the disadvantage of relying on
periodic Greens functions.

Instead, we use the following representation for the physical flow velocity:
\be
\uu = {\cal S}\nr_{\gamma}\mbf{f} + \uu_\tbox{resp},
\qquad\mbox{ where }
{\cal S}\nr_{\gamma}\mbf{f} :=
\sum_{|n|\le 1} {\cal S}_{\gamma+n\dd}\mbf{f}
~.
\label{utotrep}
\ee
The imposed flow (the first term)
involves only the vesicle and its
immediate neighbor images, as with \eqref{Dnr}.
We define the associated imposed pressure
similarly: ${\cal Q}\nr_{\gamma}\mbf{f} :=
\sum_{|n|\le 1} {\cal Q}_{\gamma+n\dd}\mbf{f}$.

The response $\uu_\tbox{resp} = \uu_\tbox{resp}[\mbf{f},\pdrive]$
%
%
is then the solution to the single-unit-cell BVP
\eqref{boxPDE}--\eqref{boxgT} with the following data
involving traces of the imposed flow on the walls:
\bea
\vv_U &=& -{\cal S}\nr_{\gamma}\mbf{f}|_U
\label{vesU}
\\
\vv_D &=& -{\cal S}\nr_{\gamma}\mbf{f}|_D
\\
\mbf{g}_u & = &
- {\cal S}_{\gamma-\dd}\mbf{f}|_R +
{\cal S}_{\gamma+\dd}\mbf{f}|_L
\\
\mbf{g}_T &=&
- T({\cal S}_{\gamma-\dd}\mbf{f}, {\cal Q}_{\gamma-\dd}\mbf{f})|_R +
T({\cal S}_{\gamma+\dd}\mbf{f}, {\cal Q}_{\gamma+\dd}\mbf{f})|_L
   + \pdrive\nn
\label{vesgT}
\eea
%
It is simple to check that \eqref{utotrep} then satisfies
no-slip velocity data on $U$ and $D$, is periodic,
and the
pressure representation has the required pressure drop \eqref{perBVPf}.
Note that, as in the $C$ block of the previous section,
there is cancellation in the discrepancies $\mbf{g}_u$ and $\mbf{g}_T$
so that even when vesicles come close to,
or intersect, $L$ or $R$, there are no near-field terms.
Effectively, the $L$ and $R$ walls are ``invisible'' to the vesicles.

To summarize, the algorithm for solving the static periodic pipe flow problem given vesicle forces $\mbf{f}$ and the driving $\pdrive$ has three main steps:
\begin{enumerate}
\item[i) ] Evaluate the right-hand side data \eqref{vesU}--\eqref{vesgT},
ie
$$
\left[ \begin{array}{c}
v\\
g\\
\end{array} \right]
\;=\;
\left[ \begin{array}{c}
-S\nr_{U,\gamma}\\ -S\nr_{D,\gamma} \\ -S_{R,\gamma-\mbf{d}}+S_{L,\gamma+\mbf{d}}\\
-K_{R,\gamma-\mbf{d}}+K_{L,\gamma+\mbf{d}}
\end{array} \right]
\mbf{f}
\,+\,
\left[ \begin{array}{c}
0\\0\\0\\
\mbf{n} \end{array} \right]
\pdrive
~;
$$
this will be done with the FMM,
except for when the vesicle is close to the wall, in which
case a recent spectral close evaluation scheme is used
\cite{lsc2d}.

\item[ii) ] Solve the rectangular linear system \eqref{bsys} for the density $\tau$ and
coefficient vector $c$; this is done with the fast direct solver to be
described in Section \ref{s:fds}.
\item[iii) ] Evaluate $\uu_\tbox{resp}$ (being
the solution to \eqref{boxPDE}--\eqref{boxgT})
using the representation \eqref{urep};
this will again be done via the FMM to get $\uu_\tbox{resp}|_\gamma$.
It is clear that $\uu_\tbox{resp}$ is linear both in $\mbf{f}$ and 
$\pdrive$.
\end{enumerate}
This three-step procedure will become one piece of the following
evolution scheme.


\subsection{Time-stepping scheme}

So far we have only described a quasi-static solution for $\uu$
driven by $\mbf{f}$ and $\pdrive$.
To close the system, we enforce no-slip conditions on the vesicle,
$ \dot{\xx} = \uu|_\gamma$, where $\cdot = \partial/\partial t$.
Substituting \eqref{utotrep} gives the first equation
in the integro-differential system of evolution equations for the vesicle dynamics, namely
\bea
\dot{\xx} &=&
S\nr_{\gamma,\gamma}\mbf{f} + \uu_\tbox{resp}[\mbf{f},\pdrive]|_\gamma \label{evolve}
\\
0 &=& \xx_s \cdot \dot{\xx}_s \label{inext}
\eea
where $\mbf{f} = -\kappa_B\xx_{ssss} + (\sigma\xx_s)_s$. Unlike other particulate systems (e.g., drops), the interfacial tension $\sigma(s,t)$ is not known {\em a priori} and needs to be determined as part of the solution.  It serves as a Lagrange multiplier to enforce the local inextensibility constraint---the second equation \eqref{inext} in this system---that the surface divergence of the membrane velocity is zero.  
This system is driven by $\pdrive$, which could vary in time (we take it as constant in our experiments).

The governing equations (\ref{evolve})-(\ref{inext}) are numerically stiff owing to the presence of high-order spatial derivatives in the bending force. As shown in \cite{ves2d}, explicit time-stepping schemes, such as the forward Euler method, suffer from a third-order constraint on the time step size, rendering them prohibitively expensive for simulating vesicle suspensions. Therefore, we use the semi-implicit scheme formulated in \cite{ves2d} with a few modifications to improve the overall numerical accuracy and stability. Given a time step size $\Delta t$ and the membrane position and tension at the $k$th time step, $(\xx^k$, $\sigma^k)$, we evolve to $({\bf x}^{k+1},\sigma^{k+1})$ by using a first-order semi-implicit time-stepping scheme on  (\ref{evolve})-(\ref{inext}). For implementational convenience, however, we treat the discretized membrane velocity, denoted with a slight abuse of notation by $\um = (\xx^{k+1} - \xx^{k})/\Delta t$, as the unknown instead of $\xx^{k+1}$. The scheme, then, is given by
\bea
\label{eq:timestep1}
\um - \, S\nr_{\gamma,\gamma}[- \Delta t\kappa_B \um_{ssss} + (\sigma^{k+1}\xx^k_s)_s] &=& S\nr_{\gamma,\gamma}[- \kappa_B \xx^k_{ssss}] + \,\uu_\tbox{resp}[- \kappa_B\xx^k_{ssss} + (\sigma^k\xx^k_s)_s,\pdrive]|_\gamma\\
\label{eq:timestep2}
\xx^k_s \cdot \um_s &=& 0. \label{inext_h}
\eea
Since the bending force is a nonlinear function of the membrane position, the standard principle of semi-implicit schemes---to treat the terms with highest-order spatial derivatives implicitly \cite{ascher-petzold98}---has been applied to the particular linearization. The tension is treated implicitly and the vesicle-channel interaction, explicitly\footnote{When the vesicle is located very close to the channel, say $\mathcal{O}(h)$ away where $h$ is the lowest distance between spatial grid points on the vesicle, a semi-implicit treatment of $\uu_\tbox{resp}$ would be more efficient since the interaction force also induces numerical stiffness in this scenario. Such a scheme, however, requires non-trivial modifications to our fast direct solver of Section \ref{s:fds}. Therefore, we postponed this exercise to future work.}. The single-layer operator $S\nr_{\gamma,\gamma}[\cdot]$ as well as the differential operator $(\cdot)_s$ are constructed using $\xx^k$. In summary, we solve the following linear system for the unknowns $(\um, \sigma^{k+1})$:
\be
\left[ \begin{array}{cc}
I + \,\Delta t\kappa_B S\nr_{\gamma,\gamma} \partial_{ssss} &
-S\nr_{\gamma,\gamma} \partial_s(\xx^k_s \cdot )\\
\xx^k_s \cdot \partial_s & 0
\end{array} \right]
\left[ \begin{array}{c} \um \\ \sigma^{k+1}
\end{array} \right]
=
\left[ \begin{array}{c}
 -\kappa_B S\nr_{\gamma,\gamma} \xx^k_{ssss} + \uu_\tbox{resp}|_\gamma
\\ 0 \end{array} \right]
\label{tsteplinsys}
\ee
with given $(\xx^k, \sigma^k)$ and then update the membrane positions as $\xx^{k+1} = \xx^k + \Delta t \um$. The operators are discretized using a spectrally-accurate Nystr\"om method
(with periodic Kress corrections for the log singularity \cite[Sec.~12.3]{LIE})
for the single-layer operator and a standard periodic spectral scheme for the differentiation operators.
The resulting discrete linear system is solved via GMRES, with all distant interactions applied using the Stokes FMM (e.g., see Appendix D of \cite{ves2d}).

On the initial time step, we generally set $\sigma_0=0$. One could obtain an improved initial guess for $\sigma_0$ by setting $\sigma_0=0$ and solving (\ref{eq:timestep1}) and (\ref{eq:timestep2}) for $\sigma$ with $\Delta t=0$.  This could be thought of as the first iteration of a Picard iteration for $\sigma$. To enable long-time simulations, we incorporate three supplementary steps in our time-stepping scheme. First, we modify the constraint equation \eqref{inext_h} as\footnote{The main idea here is similar in spirit to the correction formula applied in \cite{anna04}, but in our scheme, we do not introduce any penalty parameters and also rigorously prove its convergence rate  (Appendix A).}
\begin{equation} \xx^k_s \cdot \um_s = \frac{L_0 - L_k}{\Delta t L_k}, \label{ALC}\end{equation}
where $L_k$ represents the perimeter of the vesicle at the $k$th time step. We will refer to this as the ``arc length correction (ALC).'' We prove in Appendix A that without the ALC, the perimeter of a vesicle will increase monotonically with the number of time steps (total error still scales as $\mathcal{O}(\Delta t)$). This would mean that the vesicle's reduced area can become very low when a large number of time steps are taken. Consequently, the simulated dynamics may correspond to a totally different system than what was originally intended (e.g., a tank-treading vesicle in shear flow might tumble if the reduced area is lowered enough). Executing the ALC at every time step, on the other hand, guarantees a $\mathcal{O}(\Delta t^2)$ convergence rate, but more importantly, the error is independent of the number of time steps for a fixed $\Delta t$ (see Theorem A.2). Second, we correct the error incurred in the enclosed area of the vesicle after every time step by solving a quadratic equation in one variable (see Appendix B). Finally, we reparameterize $\gamma$ at every time step so that spatial discretization points are located {\em approximately} equal arc lengths apart (see Appendix C). 

\textbf{Suspension flow.} Although we have presented the case with a single (periodized) vesicle $\gamma$, the above scheme carries over naturally to multiple vesicles. Since such extension has been described previously in other contexts (e.g., see \cite{ves2d} for free-space and  \cite{rahimian10} for constrained geometry problems) and does not modify in any way our periodization scheme, we only highlight the main steps here. First, the single-layer potential in the representation \eqref{utotrep} is replaced with a sum of such potentials over all of the vesicles. The equations for the imposed flow data on the walls \eqref{vesU}--\eqref{vesgT} are then modified accordingly. In discretizing the evolution equation for each individual vesicle, the bending force in the self-interaction term is treated semi-implicitly, similar to (\ref{eq:timestep1})--(\ref{eq:timestep2}). However, the bending forces in the vesicle-vesicle interaction terms can either be treated semi-implicitly, resulting in a dense linear system, or explicitly, resulting in a block tri-diagonal system. While the latter scheme has a marginally smaller computational cost, the former scheme has better stability properties in general since vesicle-vesicle interactions can induce stiffness into the evolution equations when they are located close to each other. In our implementation, we treat all vesicle interactions semi-implicitly.  We use a recent single-layer close evaluation scheme \cite{lsc2d} to compute the nearby vesicle interactions, whereas for distant interactions, we use the FMM.

\section{Fast direct solver for the fixed channel geometry}
\label{s:fds}
At every time step, the right-hand side in
\eqref{tsteplinsys} must be evaluated, which involves the
channel response 3-step solution given at the end of Section~\ref{s:uresp}.
However, since the channel geometry is fixed, a fast direct solver
enables the second, potentially most expensive, of these three steps to be
performed in ${\cal O}(N)$ time with a very small constant.
Recall that this step involves
solving the $2\times 2$ block integral equation system \eqref{bsys}.
Upon discretization, this becomes a rectangular
system of size $2(N+K)\times2(N+M)$ given by
\begin{equation}
 \mtwo{\mtx{A}}{\mtx{B}}{\mtx{C}}{\mtx{Q}} \vtwo{\vct{\tau}}{\vct{c}} = \vtwo{\vct{v}}{\vct{g}}.
 \label{eq:blocksystem}
 \end{equation}

This section presents a fast direct solution technique for (\ref{eq:blocksystem}).  The idea is to 
precompute for ${\cal O}(N)$ computational cost, the factors in the block matrix $2\times 2$ 
solve.  Then, the contribution from the channel geometry in the time-stepping scheme only requires a collection of 
inexpensive linear-scaling matrix vector multiplies.

This section begins by presenting the $2\times 2$ block solve.  The remainder of the section describes 
how to efficiently build and apply the block solver.  The bulk of the novelty in this work lies in the 
linear scaling technique for representating the interactions of neighboring geometries.  

\subsection{The block solve}
The solution to the rectangular system \eqref{eq:blocksystem} is given by 

\begin{align*} 
\vct{c} &=  -\mtx{S}^\dag \left(\vct{g}-\mtx{C}\mtx{A}^{-1}\vct{v}\right)\\
 \vct{\tau}& =\mtx{A}^{-1}\vct{v}-\mtx{A}^{-1}\mtx{B}\vct{c},
\end{align*}
where $\mtx{S} = \mtx{Q}-\mtx{C}\mtx{A}^{-1}\mtx{B}$, and $\mtx{S}^\dag$ denotes the 
pseudoinverse of $\mtx{S}$.  $\mtx{S}$ is often referred to 
as a Schur complement matrix.

The Schur complement and $\mtx{A}^{-1}\mtx{B}$ need only be formed once, independent of 
the number of right-hand sides (i.e. time steps).  Beyond that, the array $\mtx{A}^{-1}\vct{v}$ need only 
be computed once per solve.

When a large number of points $N$ are needed to discretize the walls due to a complex geometry 
or small vesicle size, the cost of computing the block solve is dominated by the cost of computing 
the inverse of $\mtx{A}$.  When $N$ is large, computing $\mtx{A}^{-1}$ is computationally
prohibitive.  Fortunately, the matrix $\mtx{A}$ has structure which can be exploited to reduce the 
computational cost of the block solve. 

Recall that $\mtx{A}$ is the discretization of (\ref{Dnr}) added to $-\frac{1}{2}\mtx{I}$.  The underlying structure of $\mtx{A}$ 
can exploited first by considering its expanded form given by
\begin{equation}\label{eq:formA}\mtx{A} = \mtx{A}_0+(\mtx{A}_{-1}+\mtx{A}_1),\end{equation}
where $\mtx{A}_j$ corresponds to the discretization of the self ($j=0$) and neighbor ($j=-1,1$) channel geometry interactions.  


Since the majority of the discretization points on the neighboring geometries are well-separated,
potential theory states that $\mtx{A}_{-1}$ and $\mtx{A}_1$ are low rank (i.e. $\mtx{A}_j$ has a rank $l$, where $l \ll N$).  
Thus, each matrix admits a factorization of the form $\mtx{A}_j = \mtx{L}_j\mtx{R}_j$,
where $\mtx{L}_j$ and $\mtx{R}_j$ are of size $N\times l$ for $j = -1,1$.  Thus, 
$$\mtx{A} = \mtx{A}_0 + \mtx{L}\mtx{R}$$,
where 
$\mtx{L} = \left[\mtx{L}_{-1} |\mtx{L}_1\right] $ and 
$\mtx{R} = \left[\mtx{R}_{-1}^T |\mtx{R}_1^T\right]^T$.  

A consequence of utilizing the low rank factorization is that the inverse of $\mtx{A}$ can 
be computed via the \emph{Sherman-Morrison-Woodbury formula} 
\begin{equation}\label{eq:woodbury1}\left(\mtx{A}_0+\mtx{L}\mtx{R}\right)^{-1} = 
\mtx{A}_0^{-1}-\mtx{A}_0^{-1}\mtx{L}\left(\mtx{I}+\mtx{R}\mtx{A}_0^{-1}\mtx{L}\right)^{-1}\mtx{R}\mtx{A}_0^{-1}.
\end{equation}

Note, only the inverse of $\mtx{A}_0$ and $\left(\mtx{I}+\mtx{R}\mtx{A}_0^{-1}\mtx{L}\right)$ need to 
be computed.  Section \ref{sec:HBS} describes a technique for inverting $\mtx{A}_0$ with a 
cost that scales linearly with $N$ for most wall geometries.  Other fast inversion techniques, such 
as \cite{2009_xia_superfast,2007_shiv_sheng,2010_xia,2010_borm_book,2004_borm_hackbusch,m2011_1D_survey,HODLR},
can be utilized in place of the method in Section \ref{sec:HBS}.  
The matrix $\left(\mtx{I}+\mtx{R}\mtx{A}_0^{-1}\mtx{L}\right)$ is $2l \times 2l$ in size and 
is small enough to be inverted rapidly via dense linear algebra.

\subsection{Construction of low-rank factorization for neighbor interactions}
\label{sec:low-rank}
The cost of constructing the factorizations of $\mtx{A}_{-1}$ and $\mtx{A}_{1}$ using general linear 
algebra techniques, such as QR, is ${\cal O}(N^2l)$ and thus would negate the key reduction in 
asymptotic complexity gained by using fast direct solvers, such as the one described in 
Section \ref{sec:HBS} to invert $\mtx{A}_0$.  
To maintain the optimal asymptotic complexity, we utilize ideas from potential theory, similar to those in 
\cite{qps}.  Unlike \cite{qps}, the neighboring geometries touch the geometry in the unit cell.  As a result,
a new technique for computing the factorizations is required. For simplicity of presentation, we describe 
the new low-rank factorization technique for the matrix $\mtx{A}_1 = \mtx{L}_1\mtx{R}_1$.  
The technique is applied in a similar manner to construct the low-rank factorization 
of $\mtx{A}_{-1} = \mtx{L}_{-1}\mtx{R}_{-1}$.  As in \cite{qps}, an interpolatory decomposition \cite{gu1996,lowrank}
is utilized.

\begin{definition}
 The \textit{interpolatory decomposition} of an $m\times n$ matrix $\mtx{M}$ that has rank $l$ is
the factorization 
$$ \mtx{M} = \mtx{P}\mtx{M}(J(1:l),:),$$
where $J$ is a vector of integers $j_i$ with $1\leq j_i\leq m$, and $\mtx{P}$ is an $m\times l$ matrix that contains an $l \times l$ identity matrix.
Namely, $\mtx{P}(J(1:l),:) = \mtx{I}_l$.  
\end{definition}

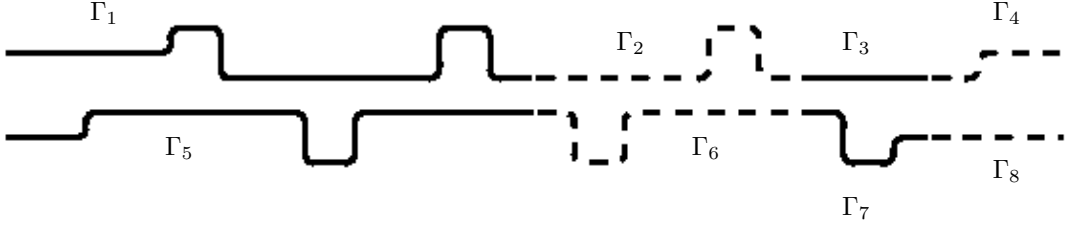
\begin{figure}
\setlength{\unitlength}{1mm}
\begin{picture}(140,50)
\put(05,00){\includegraphics[height=50mm]{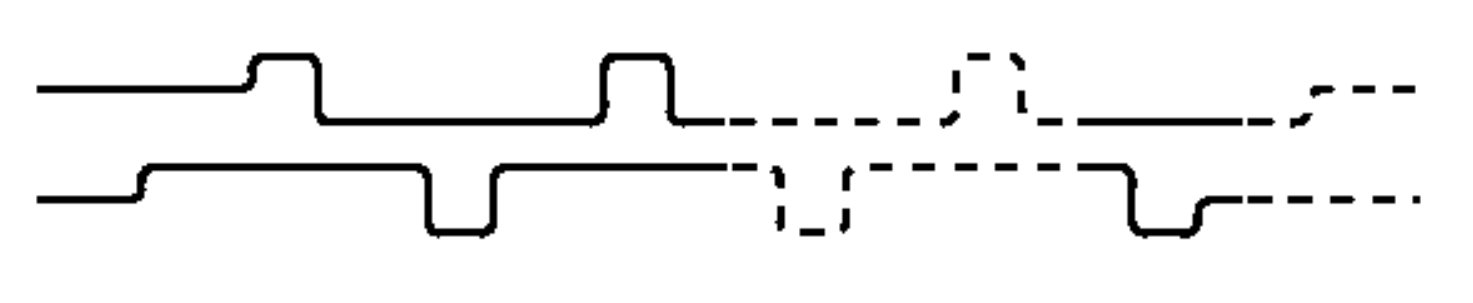}}
\put(20,26){$\Gamma_1$ }
\put(90,22){$\Gamma_2$}
\put(120,22){$\Gamma_3$}
\put(140,26){$\Gamma_4$}
\put(30,8){$\Gamma_5$ }
\put(100,8){$\Gamma_6$}
\put(120,00){$\Gamma_7$}
\put(140,5){$\Gamma_8$}
\end{picture}
\caption{Illustration of the dyadic partitioning used to compress the interaction of the walls in the unit cell
with the right neighbor, i.e. $\mtx{A}_{1}$.  }
\label{fig:dyadic}
\end{figure}

First, the upper and lower part of the geometry is partitioned into a 
collection of $M$ sections $\Gamma_j$ via dyadic refinement where the rectangular boxes 
enclosing each segment get smaller as they approach the neighboring cells. Thus, $\displaystyle\partial \Omega_0 =\cup_{j=1}^M \Gamma_j$,
where $M$ is the number of sections $\partial \Omega_0$ is partitioned into.
Figure \ref{fig:dyadic} illustrates the partitioning of the walls for the channel with reservoirs 
in Figure \ref{fig:geom}(d) when compressing the interaction with the right neighbor $\Omega_1$.  
The refinement is stopped when the smallest box contains no more than a specified number 
of points $n_{\rm max}$.  Typically, $n_{\rm max} = 45$ is a good choice.

\begin{figure}[t]
\centering
	\includegraphics[width=\textwidth]{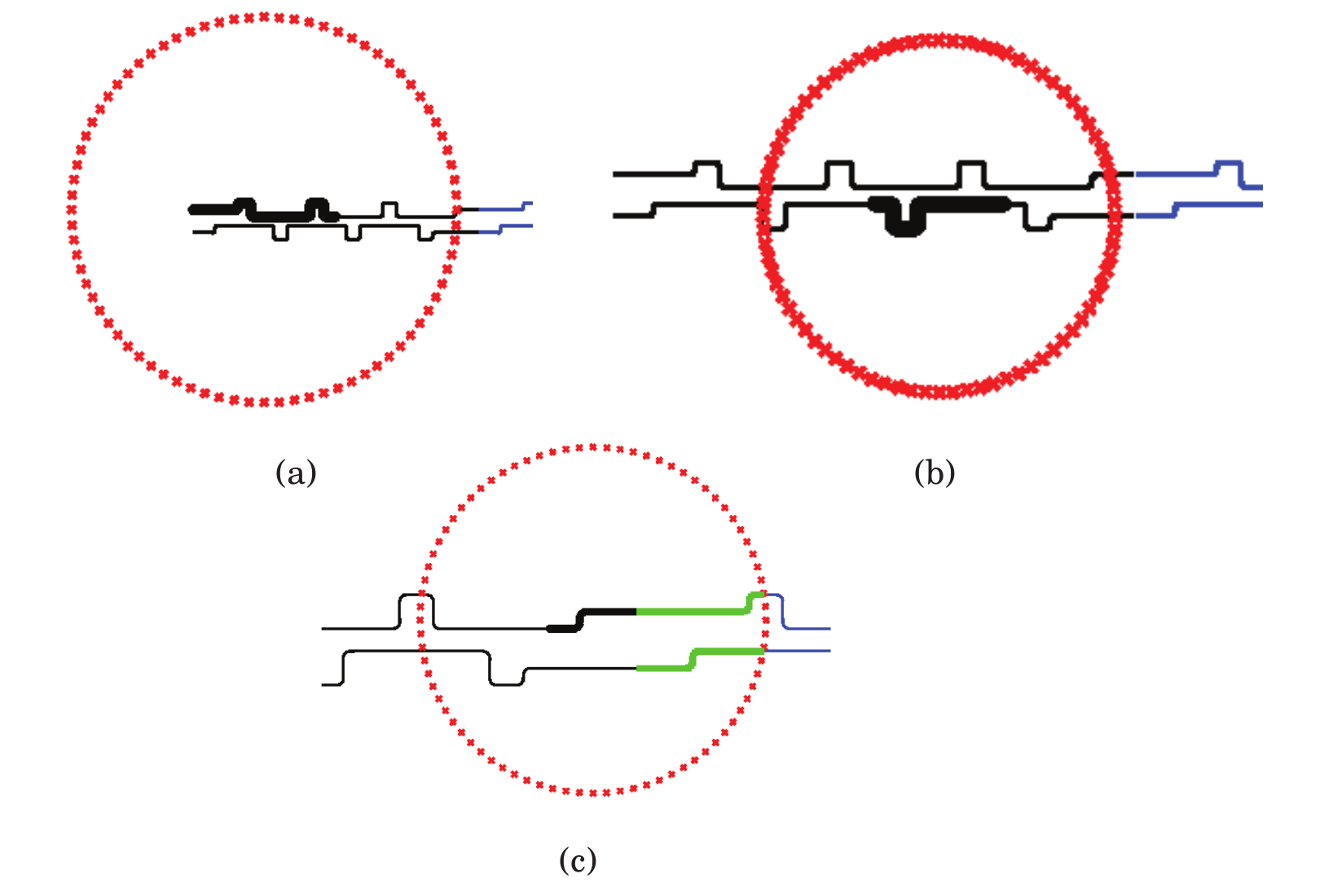}%
\caption{(a) and (b) illustration of the proxy points ({\color{red} \textbf{x}}) for $\Gamma_1$, $\Gamma_6$.  (c) illustrates
the proxy points and the near points (thick green line) on $\partial \Omega_1$ for $\Gamma_4$.  }
\label{fig:proxy}
\end{figure}

For each portion of the boundary $\Gamma_j$ that is not touching $\Omega_1$, consider a circle concentric 
with a box bounding $\Gamma_j$
with radius slightly less than the distance from the center of the box 
enclosing $\Gamma_j$ to the ``wall'' where $\Omega_0$ and $\Omega_1$ meet.  
From potential theory, we know that any field generated by sources outside 
of this circle can be approximated well by placing enough equivalent 
charges on the circle.  In practice, it is possible to place a small number 
of ``proxy'' points spaced evenly on the circle.  For the experiments
 in this paper, we found it is sufficient to have $80$ proxy points.  
 Figure \ref{fig:proxy} (a) and (b) illustrate the proxy points for 
 $\Gamma_1$ and $\Gamma_6$.  
An interpolatory decomposition is then found for the matrix denoted $A^{\rm proxy}$ 
that characterizes the interactions between $\Gamma_j$ and the proxy 
points.  The result is a matrix $P_j$ and 
index vector $J_j$.  For $\Gamma_j$ touching 
$\partial\Omega_1$, a collection of $80$ proxy points are placed evenly on a
circle with a radius 1.75 times the radius of the smallest circle enclosing $\Gamma_j$.  
All of the points on $\partial\Omega_1$ that lie inside this proxy circle
are called \emph{near} points.  Figures \ref{fig:proxy}(c) and (d) illustrates the proxy and 
near points for $\Gamma_4$.  An interpolatory decomposition is then formed for the matrix 
$$[\mtx{A}_1(\Gamma_j,I^{\rm near}) | A^{\rm proxy}]$$, where 
$I^{\rm near}$ corresponds to the indices of the portion of $\partial\Omega_1$ 
that are near $\Gamma_j$.  The points $\partial\Omega_0$ picked by the 
interpolatory decomposition are called \emph{skeleton points}.


We can now assemble the factors $\mtx{L}_1\mtx{R}_1$ by sweeping through 
the regions $\Gamma_j$. Let $J= [J_1(1:l_1),\ldots,J_M(1:l_M)]$.  
Then the matrix $\mtx{L}_1$ is a block diagonal matrix where each block is a $\mtx{P}_j$ 
and $\mtx{R}_1 = \mtx{A}_1(J,:)$.

\begin{remark}
 The factorization as described does not have optimal rank.  Optimal rank 
 can be obtained by recompressing via additional low rank factorizations 
 while assembling $\mtx{L}_1$ and $\mtx{R}_1$.  Depending on the geometry and the
 computer, it may or may not be beneficial to do the recompression.
\end{remark}

\subsection{HBS inversion}
\label{sec:HBS}
As stated previously, the dense matrix $\mtx{A}_0$ has structure that we call 
\textit{Hierarchically Block Separable (HBS)}.  The HBS structure allows for 
an approximation of $\mtx{A}_0^{-1}$ to be computed rapidly.
Loosely speaking, its off-diagonal blocks are low rank.
This arises because $\mtx{A}_0$ is the discretization on a curve of an
integral operator with smooth kernel (when the walls are not space-filling).  
This section briefly describes the HBS property and how it can be exploited 
to rapidly construct an approximate inverse 
of a matrix. For additional details, see \cite{m2011_1D_survey}.
Note that the HBS property is very similar to the concept of \textit{Hierarchically
Semi-Separable (HSS)} matrices \cite{2007_shiv_sheng,gu_divide}.

\subsubsection{Block separable}
Let $\mtx{M}$ be an $mp\times mp$ matrix that is blocked into $p\times p$ blocks,
each of size $m\times m$.

We say that $\mtx{M}$ is ``block separable'' with ``block-rank'' $k$
if for $\tau = 1,\,2,\,\dots,\,p$, there exist $n\times k$
matrices $\mtx{U}_{\tau}$ and $\mtx{V}_{\tau}$ such that each off-diagonal
block $\mtx{M}_{\sigma,\tau}$ of $\mtx{M}$ admits the factorization
\begin{equation}
\label{eq:yy1}
\begin{array}{cccccccc}
\mtx{M}_{\sigma,\tau}  & = & \mtx{U}_{\sigma}   & \tilde{\mtx{M}}_{\sigma,\tau}  & \mtx{V}_{\tau}^{*}, &
\quad \sigma,\tau \in \{1,\,2,\,\dots,\,p\},\quad \sigma \neq \tau.\\
m\times m &   & m\times k & k \times k & k\times m
\end{array}
\end{equation}

Observe that the columns of $\mtx{U}_{\sigma}$ must form a basis for
the columns of all off-diagonal blocks in row $\sigma$, and
analogously, the columns of $\mtx{V}_{\tau}$ must form a basis for the
rows in all of the off-diagonal blocks in column $\tau$. When (\ref{eq:yy1})
holds, the matrix $\mtx{M}$ admits a block factorization

\begin{equation}
\label{eq:yy2}
\begin{array}{cccccccccc}
 \mtx{M} &  =& \mtx{U}&\tilde{\mtx{M}}& \mtx{V}^{*} & +&  \mtx{D},\\
mp\times mp &   & mp\times kp & kp \times kp & kp\times mp && mp \times mp\\
\end{array}
\end{equation}

where
$$
\mtx{U} = \mbox{diag}(\mtx{U}_{1},\,\mtx{U}_{2},\,\dots,\,\mtx{U}_{p}),\quad
\mtx{V} = \mbox{diag}(\mtx{V}_{1},\,\mtx{V}_{2},\,\dots,\,\mtx{V}_{p}),\quad
\mtx{D} = \mbox{diag}(\mtx{D}_{1},\,\mtx{D}_{2},\,\dots,\,\mtx{D}_{p}),
$$

and

$$\tilde{\mtx{M}} = \left[\begin{array}{cccc}
0 & \tilde{\mtx{M}}_{12} & \tilde{\mtx{M}}_{13} & \cdots \\
\tilde{\mtx{M}}_{21} & 0 & \tilde{\mtx{M}}_{23} & \cdots \\
\tilde{\mtx{M}}_{31} & \tilde{\mtx{M}}_{32} & 0 & \cdots \\
\vdots & \vdots & \vdots
\end{array}\right].
$$

Once the matrix $\mtx{M}$ has been put into block separable form, its inverse is given by 
\begin{equation}
\label{eq:woodbury}
\mtx{M}^{-1} = \mtx{E}\,(\tilde{\mtx{M}} + \hat{\mtx{D}})^{-1}\,\mtx{F}^{*} + \mtx{G},
\end{equation}
where
\begin{align}
\label{eq:def_muhD}
\hat{\mtx{D}} =&\ \bigl(\mtx{V}^{*}\,\mtx{D}^{-1}\,\mtx{U}\bigr)^{-1},\\
\label{eq:def_muE}
\mtx{E}  =&\ \mtx{D}^{-1}\,\mtx{U}\,\hat{\mtx{D}},\\
\label{eq:def_muF}
\mtx{F}  =&\ (\hat{\mtx{D}}\,\mtx{V}^{*}\,\mtx{D}^{-1})^{*},\\
\label{eq:def_muG}
\mtx{G}  =&\ \mtx{D}^{-1} - \mtx{D}^{-1}\,\mtx{U}\,\hat{\mtx{D}}\,\mtx{V}^{*}\,\mtx{D}^{-1}.
\end{align}

\subsubsection{Hierarchically Block-Separable}
Informally speaking, a matrix $\mtx{M}$ is \textit{Hierarchically Block-Separable} (HBS)
if it is amenable to a \textit{telescoping} version of the above
block factorization.  In other words,
in addition to the matrix $\mtx{M}$ being block separable, so is $\tilde{\mtx{M}}$
once it has been re-blocked to form a matrix with $p/2 \times p/2$ blocks,
and one is able to repeat the process in this fashion multiple times.


For example, a ``3 level'' factorization of $\mtx{M}$ is
\begin{equation}
\label{eq:united4}
\mtx{M} = \mtx{U}^{(3)}\bigl(\mtx{U}^{(2)}\bigl(\mtx{U}^{(1)}\,\tilde{\mtx{M}}^{(0)}\,(\mtx{V}^{(1)})^{*} + \mtx{B}^{(1)}\bigr)
(\mtx{V}^{(2)})^{*} + \mtx{B}^{(2)}\bigr)(\mtx{V}^{(3)})^{*} + \mtx{D}^{(3)},
\end{equation}
where the superscript denotes the level.

The HBS representation of an $N\times N$ matrix requires $\mathcal{O}(Nk)$ to store and to apply to a vector.
By recursively applying formula (\ref{eq:woodbury}) to the telescoping factorization, an approximation of the inverse can be computed with $\mathcal{O}(Nk^2)$ computational cost; see \cite{m2011_1D_survey}.
This compressed inverse can be applied to a vector (or a matrix) very rapidly.

\section{Numerical results}
\label{s:results}
\begin{algorithm}[!t]
\caption{Main (single vesicle)}\label{Main}
\begin{algorithmic}
\State \emph{Step 1:} Compress $A^{-1}$ using the fast direct solver
\For{$k=0 : N_{\Delta t} - 1$}
\State \emph{Step 2:} Compute bending and tension forces (${\bf f}_\sigma={\bf 0}$ if $k=0$)
\begin{equation}
{\bf f}_b=-\kappa_B {\bf x}_{ssss}^k,\quad {\bf f}_\sigma=\left(\sigma^k{\bf x}_s^k\right)_s,\quad {\bf f}={\bf f}_b+{\bf f}_\sigma
\end{equation}
\State \emph{Step 3:} Compute vesicle-wall and vesicle-side interactions (step (i) from Section \ref{s:uresp})
\begin{align}
{\bf v}_U&=\left.-\mathcal{S}_\gamma^\text{near}{\bf f}\right|_U\\
{\bf v}_D&=\left.-\mathcal{S}_\gamma^\text{near}{\bf f}\right|_D\\
{\bf g}_u&=\left.-\mathcal{S}_{\gamma-{\bf d}}{\bf f}\right|_R+\left.\mathcal{S}_{\gamma+{\bf d}}{\bf f}\right|_L\\
{\bf g}_T&=\left.-T\left(\mathcal{S}_{\gamma-{\bf d}}{\bf f},\mathcal{Q}_{\gamma-{\bf d}}{\bf f}\right)\right|_R+\left.T\left(\mathcal{S}_{\gamma+{\bf d}}{\bf f},\mathcal{Q}_{\gamma+{\bf d}}{\bf f}\right)\right|_L+p_\text{drive}{\bf n}
\end{align}
\State \emph{Step 4:} Solve for $\bm \tau$, ${\bf c}$ using the fast direct solver (step (ii) from Section \ref{s:uresp})
\State \emph{Step 5:} Compute wall-vesicle and proxy-vesicle interactions (step (iii) from Section \ref{s:uresp})
\begin{equation}
{\bf u}_\text{resp}=\left.\mathcal{D}_\Gamma^\text{near}{\bm \tau}\right|_\gamma+\sum_{m=1}^M{\left.{\bf c}_m\phi_m\right|_\gamma}
\end{equation}
\State \emph{Step 6:} Solve (\ref{tsteplinsys}) for ${\bf u}$ and $\sigma^{k+1}$ using GMRES with the constraint equation modified as \eqref{ALC}. 
\State \emph{Step 7:} Set ${\bf x}^{k+1}={\bf x}^k+\Delta t\:{\bf u}$
\State \emph{Step 8:} Apply area correction to ${\bf x}^{k+1}$ (Appendix B)
\State \emph{Step 9:} Apply the reparameterization scheme on ${\bf x}^{k+1}$ (Algorithm 2 of Appendix C)

\EndFor
\end{algorithmic}
\label{alg:main}
\end{algorithm}

\begin{figure}[t]
\centering
	\includegraphics[width=\textwidth]{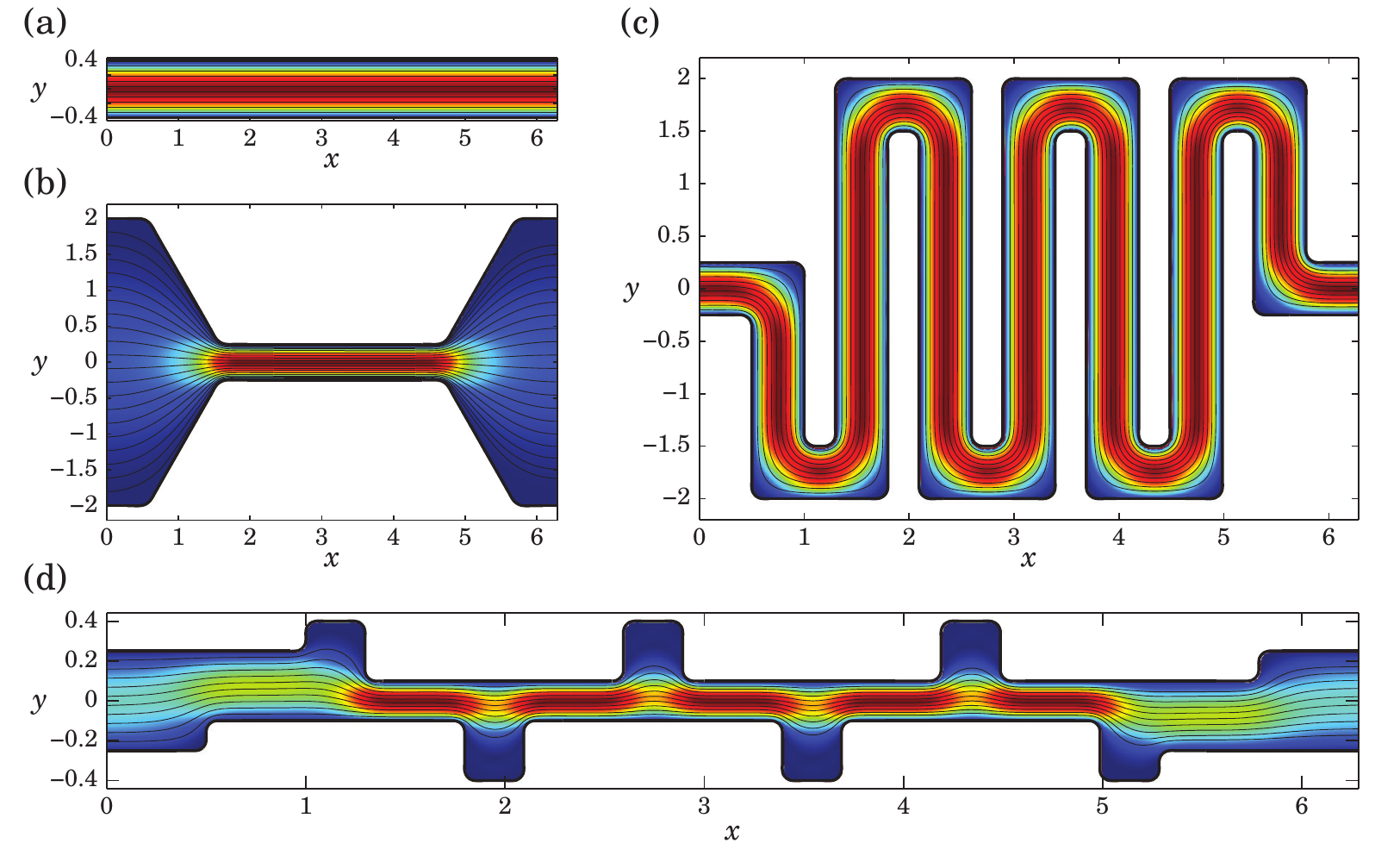} 
 	\caption{Streamlines for pressure-driven flows through four geometries with no-slip boundary conditions. The background color indicates the magnitude of the velocity (red indicates high and blue indicates low). For the remainder of this paper, we will refer to geometry (a) as the straight channel, (b) as the converging-diverging channel, (c) as the serpentine channel, and (d) as the channel with reservoirs. In all four geometries, we use $N=4{,}000$ quadrature points on each wall, $K=32$ points on each side, and $M=128$ proxy points. The ring of proxy sources has radius $d=2\pi$.}
\label{fig:geom}
\end{figure}

In this section, we test the performance of our algorithm on simulating flows through four geometries ranging in complexity from a flat channel to a more complicated space-filling geometry. We apply a pressure-drop on each side which drives the flow in the positive $x$-direction while simultaneously enforcing no-slip boundary conditions on the walls. The streamlines are plotted in Fig.~\ref{fig:geom} along with the magnitude of the velocity which is indicated by the color in the background. We will refer to geometry (a) as the straight channel, (b) as the converging-diverging channel, (c) as the serpentine channel, and (d) as the channel with reservoirs. Note that all the wall geometries are $\mathcal{C}^\infty$ curves, constructed using partitions of unity to avoid sharp corners. We used the trapezoidal rule, with $N=4{,}000$ nodes on each wall, to evaluate the layer potentials at interior targets for this test\footnote{For some applications, the trapezoidal rule may not be an ideal choice since it cannot yield uniform convergence for points that are very close to the walls. For example, notice the thin band near the wall in Fig.~\ref{fig:serperror}. Accurate close-evaluation schemes must be used instead. In our setting, the no-slip condition at the walls leads to low velocities near the walls and in general helps alleviate numerical instabilities naturally. Incorporating close evaluation schemes (for walls) is one of the first tasks in our future work.}\\
\indent
In the following subsection, we analyze the convergence of the periodization scheme as we vary the number of proxy points (denoted by $M$), discretization points on $L$ and $R$ (denoted by $K$), and quadrature points on the walls (denoted by $N$). For these tests, we set the boundary conditions and pressure-drop such that the horizontal Poiseuille flow ${\bf u}_e({\bf x})=(\alpha x_2^2,0)$, $p_e({\bf x})=\alpha \mu x_1$, with $\alpha=p_\text{drive}/(2\mu d)$ satisfies the BVP. For all tests, we use $\alpha=0.2$, $\mu=0.7$, and $d=2\pi$. The relative errors are defined as
\begin{align}
\text{relative error in velocity}&=\frac{\|{\bf u}-{\bf u}_e\|_2}{\max\limits_{{\bf x}\in\Omega}{\left(\| {\bf u}_e\|_2\right)}}\\
\text{relative error in pressure}&=\left|\frac{p-p_e}{p_\text{drive}}\right|,
\end{align}
where $({\bf u},p)$ is the numerical solution obtained by the algorithm. We report the maximum relative error at the same three points $(0.8,0)$, $(3.2,0)$, and $(4.8,0)$ for all four geometries (not too close to the walls).

 \begin{figure}[!b]
\centering
\includegraphics[width=\textwidth]{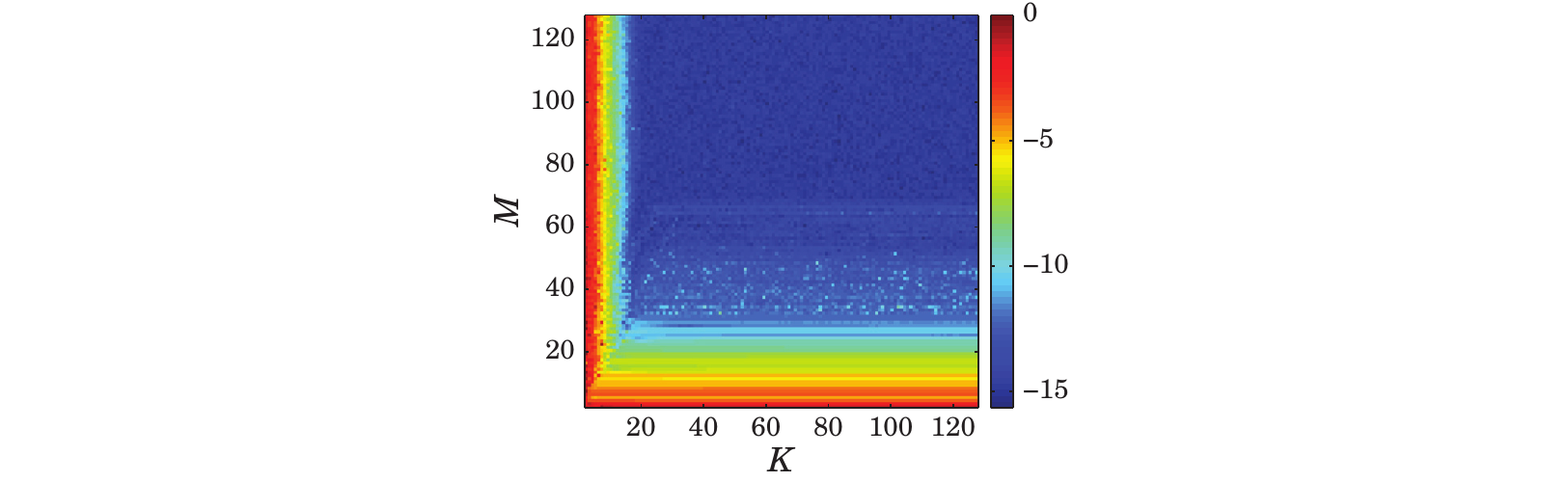} 
 \caption{Logarithm of the relative errors for the velocity as we vary the number of points on the sides and on the ring of proxy sources for the converging-diverging geometry in Fig.~\ref{fig:geom}(b). For each test point, the velocity was computed using $N=600$ quadrature points per wall. The boundary integral equation was solved using GMRES.}
 \label{fig:proxyerror}
\end{figure}

 \subsection{Periodization scheme}
Here, we analyze the convergence properties of the periodization scheme. We begin by focusing on the sides $L$ and $R$ and the proxy sources $\mathcal{C}$. We compute the relative errors for the velocity in the converging-diverging channel as the number of side points $K$ and proxy sources $M$ are varied. The number of quadrature points along the walls is fixed to $N=600$, which, based on previous experimentation, is enough to resolve the flow to near machine precision. As can be seen in Fig.~\ref{fig:proxyerror}, rapid convergence is achieved and high accuracy is obtained using a relatively small $K$ and $M$. This is mainly because of the analytical cancellation of close interactions of $L$, $R$, and $\mathcal{C}$ with the walls, as shown in equations (\ref{BIEgu}) and (\ref{BIEgT}). Since flow far from a disturbance is generally smooth with rapidly decaying Fourier modes, the distant interactions relating to $L$, $R$, and $\mathcal{C}$ are easily resolved.

\begin{figure}[!t]
\centering
 \includegraphics[width=\textwidth]{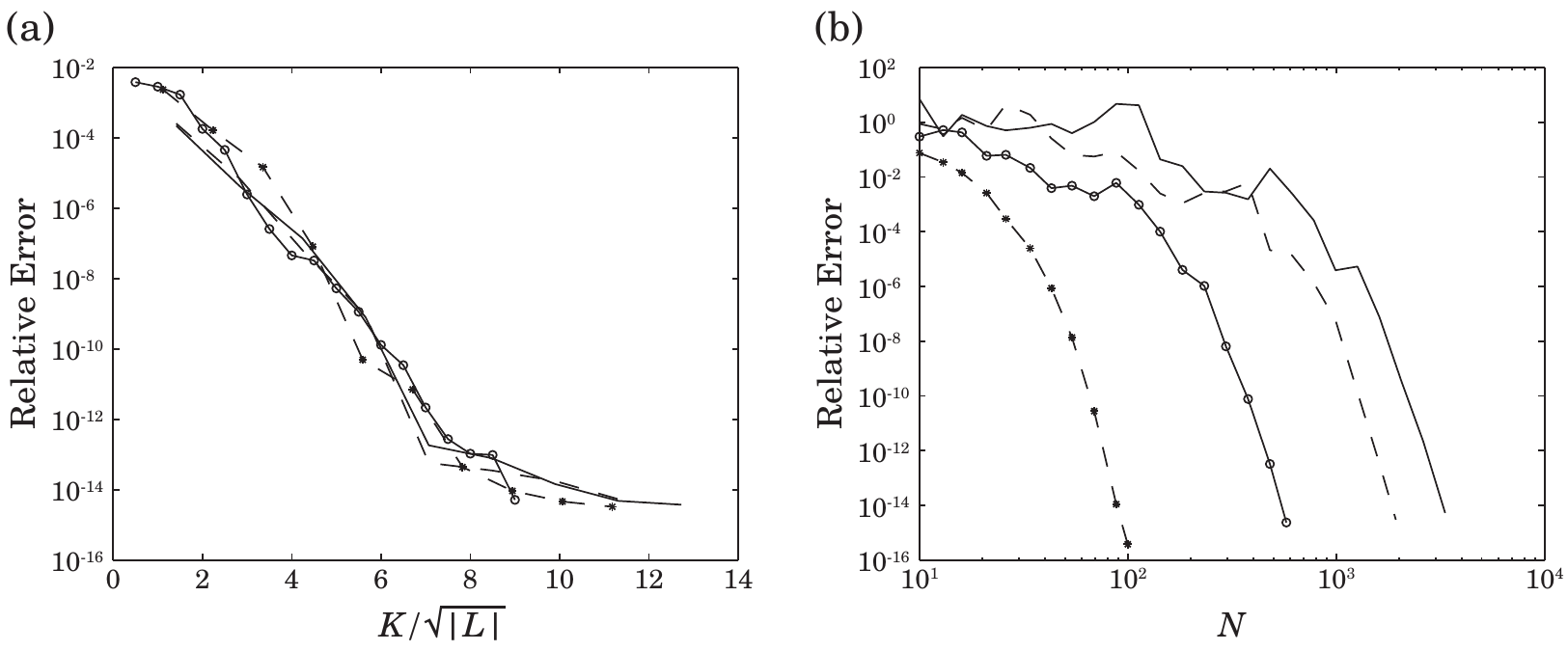} 
 \caption{$(*)$ straight channel, $(\circ)$ converging-diverging channel, $(-)$ serpentine channel, $(--)$ channel with reservoirs. (a) Relative errors as the number of side points and proxy sources are varied. For each geometry, we set the number of points on the walls to $N=4{,}000$ and the number of proxy sources to $M=4K$. The height of the inlet/outlet is denoted by $|L|$. Notice that the convergence appears to be independent of the complexity of the channel. (b) Relative errors as the number of points on the walls are varied. In each case, we used $K=32$ points per side and $M=128$ proxy sources. Spectral convergence was obtained for all four geometries.}
 \label{fig:converge}
\end{figure}

Next, we investigate how the complexity of the geometry affects the convergence rate of the side points and the proxy sources. This is done by setting $M=4K$ and measuring the relative errors pertaining to the four test geometries as $K$ is varied. Fig.~\ref{fig:converge}(a) shows the results. Notice that once we account for the height of the inlet/outlet, the convergence of all four geometries is almost exactly the same, independent of the complexity of the channel. This feature is especially nice for problems dealing with multi-scale physics, such as our application to vesicle flows.

\begin{figure}[!t]
\centering
           \includegraphics[width=\textwidth]{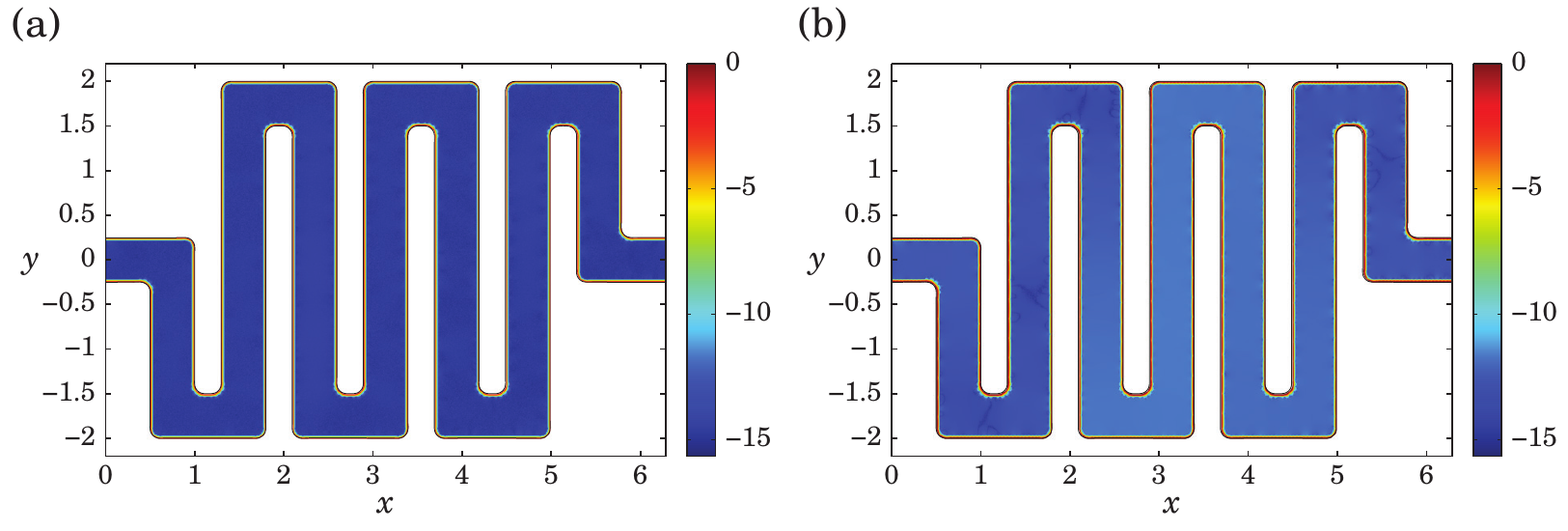} 
 \caption{Spatial plot of the logarithm of relative errors in velocity (a) and pressure (b). Errors were computed using the trapezoidal rule with 4{,}000 quadrature points per wall. The thin band near the walls, where low accuracy is obtained, may be resolved using a close evaluation scheme.}
\label{fig:serperror}
\end{figure}

Finally, we focus on the walls $U$ and $D$. Fig.~\ref{fig:converge}(b) shows the relative errors for all four test geometries as the number of quadrature points $N$ is varied. In all cases, we observe spectral convergence. (Recall that the geometries from Fig.~{\ref{fig:geom}} were constructed using partitions of unity and are $C^\infty$.) The relative errors for both the velocity and pressure inside the serpentine channel are plotted in Fig.~\ref{fig:serperror}. The small band of errors near the wall is the result of using the trapezoidal rule as opposed to a close evaluation scheme for the walls. It will turn out that in some cases this error band may cause problems for the vesicle simulation if vesicles drift too close to the walls, as we will discuss in the following section.

\begin{figure}[!t]
\centering
           \includegraphics[width=\textwidth]{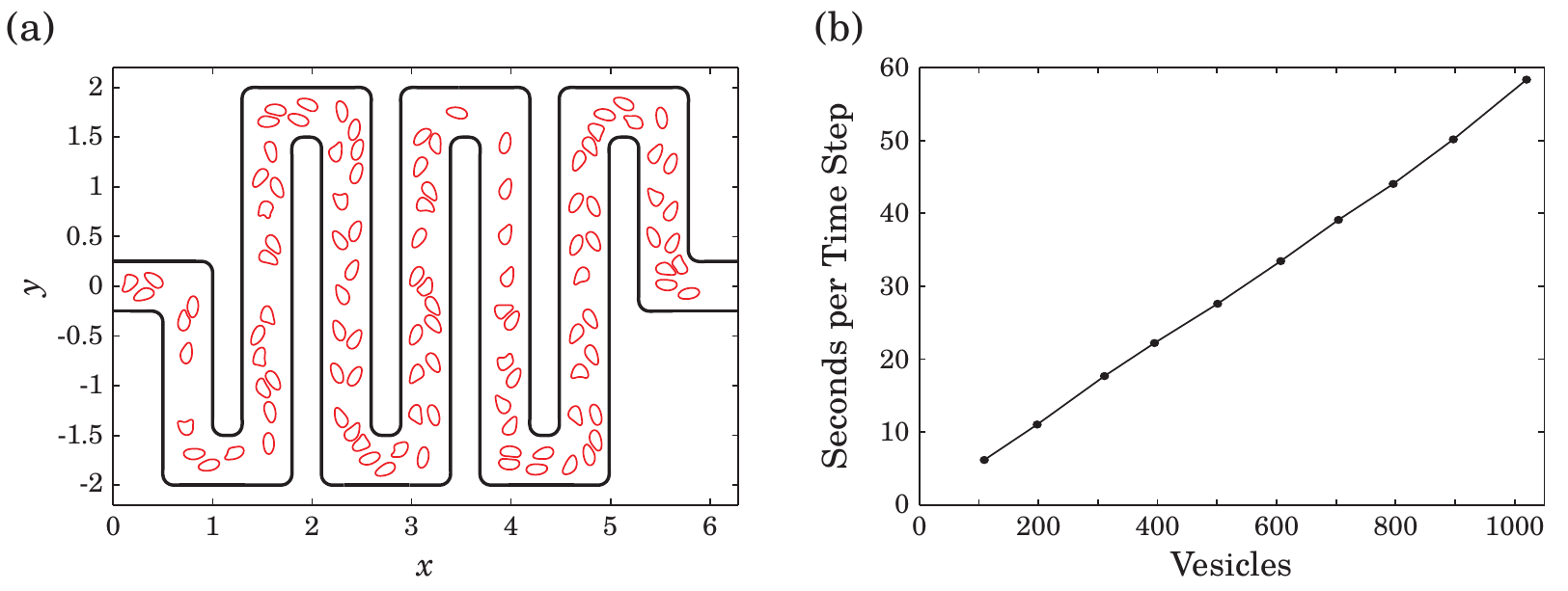} 
 \caption{(a) A snapshot of 109 vesicles in the serpentine channel. A pressure difference is pushing the vesicles in the positive $x$-direction. (b) Average time per time step (in seconds) for the first 10 time steps as the number of vesicles is varied. Each dot represents a data point. Timings were performed using a laptop with a 2.4 GHz dual-core Intel Core i5 processor and 8 GB of RAM.}
\label{fig:serpves}
\end{figure}

\begin{table}[!b]
\caption{The CPU time distribution for the first 10 time steps of a simulation with 1{,}020 vesicles (64 points each) in the serpentine channel with $N=29{,}580$ points per wall. Each time step took an average of 58 seconds on a laptop with a 2.4 GHz dual-core Intel Core i5 processor and 8 GB of RAM.}
\centering
{\setlength{\extrarowheight}{1pt}
\begin{tabular}{lc} \hline \hline
Operation & Percentage\\
\hline
vesicle to vesicle interactions & 63.49\\
vesicle to wall interactions & 14.77\\
preconditioner & 7.30\\
wall to vesicle interactions & 3.80\\
bending and tension forces & 3.17\\
inextensibility operator & 1.41\\
solve for $\tau$, $c$ using $A^{-1}$ & 1.40\\
vesicle to side interactions & 1.31\\
proxy point to vesicle interactions & 1.23\\
vesicle area corrections & 1.16\\
other & 0.96 \\ \hline \hline
\end{tabular}}
\label{tab:CPUtime}
\end{table}

\subsection{Vesicle flow simulation}
We now give numerical results for the vesicle flow simulation described in Section \ref{s:appl}. A snapshot of 109 vesicles flowing through the serpentine channel is shown in Fig.~\ref{fig:serpves}(a), where a pressure difference between the left and right sides drives the flow in the positive $x$-direction. When modeling such flows, it's vitally important that vesicles avoid collisions with other vesicles and the walls. We have found that using 64 points per vesicle along with a close evaluation scheme is usually enough to prevent vesicle-vesicle collisions as long as vesicle shapes are not elongated and the time step is not too big. For the simulations in Figs.~\ref{fig:largesim}(a) and \ref{fig:serpves}(a), the time step was set to $\Delta t=0.005$. When modeling pressure driven flows with no-slip boundary conditions, the vesicles often keep a safe distance from the walls and a close evaluation scheme is not always required. However, there are plenty of situations where this is not the case. In general, vesicle-wall collisions tend to occur more often when the channel geometry has bends or tight spaces, the number of vesicles is large, or the bending modulus $\kappa_B$ is high.

We now focus on the performance of the algorithm. Fig.~\ref{fig:serpves}(b) shows the scaling as the number of vesicles in the serpentine geometry is increased. When performing the timings, the vesicle dimensions were scaled to maintain the same relative spacing. Each vesicle had 64 discretization points and the number of quadrature points on the walls was set to 29 times the number of vesicles. The algorithm maintained linear scaling to over 1{,}000 vesicles on a laptop with a 2.4 GHz dual-core Intel Core i5 processor and 8 GB of RAM. In Table \ref{tab:CPUtime}, we give the CPU time distribution for the simulation with 1{,}020 vesicles. Approximately $82\%$ of the computational time was spent computing vesicle-vesicle, vesicle-wall, or wall-vesicle interactions. The majority of this expense was handled by the FMM. Approximately $7\%$ of the time was used to construct a preconditioner which was a sparse block diagonal matrix consisting of vesicle-vesicle self interactions.
Notice that while almost one half of the points 
are associated with the fixed geometry,
the solve associated with these points takes less than 2\% of the total time
for the time step, thanks to the fast direct solver.

\begin{figure}[t!]
\centering
\includegraphics[width=\textwidth]{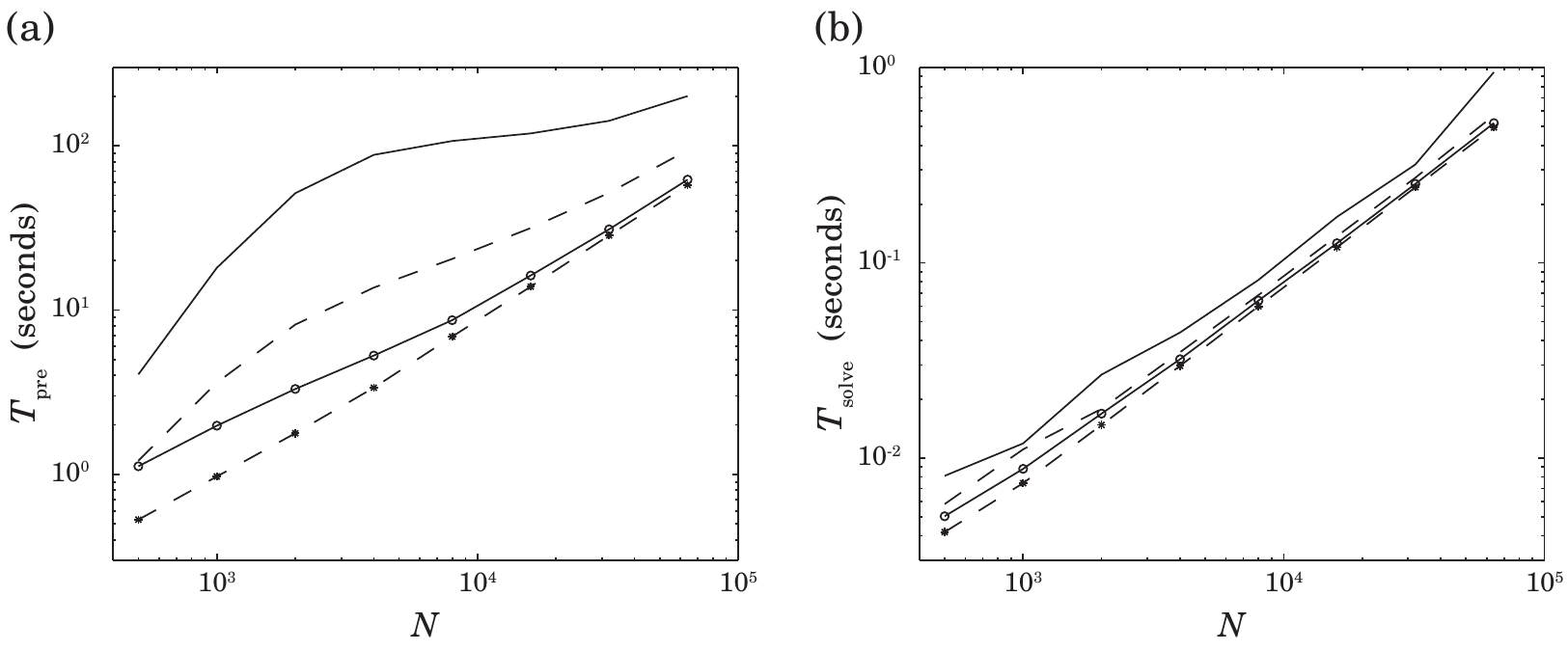} 
\caption{$(*)$ straight channel, $(\circ)$ converging-diverging channel, $(-)$ serpentine channel, $(--)$ channel with reservoirs. Time in seconds for the precomputation (a) and solve (b) steps verses the number of discretization points $N$ on the walls when the fast direct solver is applied to the geometries in Fig. \ref{fig:geom}.}
\label{fig:scaling}
\end{figure}

 \subsection{Fast direct solver}
This section reports on the performance of the fast direct solver for the wall computations.  The experiments in this 
section were performed on a laptop computer with two quad core Intel Core i7-4700MQ processors and 16 GB of RAM.  The 
user prescribed tolerance was set to $\epsilon = 10^{-10}$.  

Table \ref{tab:scaling} reports the time for the precomputation
$T_{\rm pre}$, the time for a solve $T_{\rm solve}$, and the absolute error $E$ for a point in the interior of the channel.
For all geometries, the fast direct solver scales linearly.  As the serpentine
channel is under-discretized until $N = 4{,}000$, the precomputation step in the solver does not observe asymptotic 
complexity until $N = 32{,}000$.  However, the solve step of the solver observes asymptotic complexity earlier.  

Figure~\ref{fig:scaling} reports the time in seconds for the precomputation (a) and the solve (b) steps verses the number of discretization points $N$ on the walls for all of the geometries.  The time for the precomputation of the straight channel  and the channel with reservoirs is about the same since these geometries have approximately the same rank interactions.  The rank interactions for the converging-diverging channel are smaller thus the constant for the precomputation is smaller.  Thanks to the small constant associated with applying an HBS matrix, the cost for the solves is approximately the same independent of geometry.  

In Table \ref{tab:fastcompare}, we perform a comparison between the LU decomposition, GMRES, and the fast direct solver when computing $\bm \tau$ and ${\bf c}$ for the serpentine channel. In this example, we set $\epsilon=10^{-12}$ for the fast direct solver. The timings (in seconds) do not include the initial LU factorization or the fast direct solver's compression of $A^{-1}$. Only the solve times are reported. We observe similar timings between the LU factorization and the fast direct solver for a relatively small number of points ($N\approx 400$). For larger values, the fast solver provides superior performance. In the comparison with GMRES, we used the FMM to apply $A$ and reported the timings for the system without preconditioning. We found that GMRES required approximately 259 iterations to obtain a relative error of $10^{-12}$. The number of iterations was independent of $N$ and mainly depended on the complexity of the geometry. Notice that even a single GMRES iteration requires more time than the fast solver.

\begin{table}[!t]
\caption{Time in seconds for precomputation $T_{\rm pre}$ and solve $T_{\rm solve}$ as the number of discretization points $N$ increases on the walls for the four geometries.  The absolute error $E$ at a point in the channel is also reported. }
\centering
{\setlength{\extrarowheight}{2pt}
 \begin{tabular}{l|ccc|ccc}
 \hline \hline
 &   \multicolumn{3}{c|}{Straight Channel} &  \multicolumn{3}{c}{Converging-Diverging Channel} \\ \hline
  $N$  &   $T_{\rm pre}$ & $T_{\rm solve}$ & $E$ &   $T_{\rm pre}$ & $T_{\rm solve}$ & $E$ \\ \hline
$500$&	$0.53$&	$4.18\times 10^{-3}$&$2.79\times 10^{-11}$     &$1.12$&	$5.04\times 10^{-3}$&$1.26\times 10^{-6}$\\
$1000$&	$0.97$&	$7.45\times 10^{-3}$&$1.10\times 10^{-9}$&	$1.98$&	$8.81\times 10^{-3}$&$1.12\times 10^{-9}$\\
$2000$&	$1.78$&	$1.48\times 10^{-2}$&	$1.05\times 10^{-9}$&	$3.31$&	$1.69\times 10^{-2}$&	$1.33\times 10^{-9}$\\
$4000$&	$3.37$&	$2.97\times 10^{-2}$&	$6.05\times 10^{-9}$&	$5.29$&	$3.21\times 10^{-2}$&	$7.28\times 10^{-9}$\\
$8000$&	$6.9$&	$5.99\times 10^{-2}$&	$1.68\times 10^{-8}$&	$8.68$&	$6.41\times 10^{-2}$&	$2.86\times 10^{-9}$\\
$16000$&$13.9$&	$1.20\times 10^{-1}$&	$1.13\times 10^{-8}$&	$16.25$&$1.26\times 10^{-1}$&	$1.00\times 10^{-8}$\\
$32000$&$28.5$&	$2.44\times 10^{-1}$&	$7.83\times 10^{-8}$&	$31.06$&$2.54\times 10^{-1}$&	$9.66\times 10^{-9}$\\
$64000$&$57.7$&	$4.94\times 10^{-1}$&	$1.92\times 10^{-7}$&	$62.26$&$5.20\times 10^{-1}$&	$8.28\times 10^{-9}$\\ \hline \hline
\multicolumn{7}{c}{}\\

\hline \hline
 & \multicolumn{3}{c|}{Serpentine Channel} & \multicolumn{3}{c}{Channel with Reservoirs} \\ \hline
  $N$  &   $T_{\rm pre}$ & $T_{\rm solve}$ & $E$ &  $T_{\rm pre}$ & $T_{\rm solve}$ & $E$ \\ \hline
$500$&	$4.06$&	$8.12\times 10^{-3}$&$8.62\times 10^{-2}$ & $ 1.21$&	$5.82\times 10^{-3}$&	$1.53$\\
$1000$&	$18.1$&	$1.19\times 10^{-2}$&	$2.72\times 10^{-2}$ &$3.62$&	$1.11\times 10^{-2}$&	$3.20\times 10^{-3}$\\
$2000$&	$51.5$&	$2.68\times 10^{-2}$&	$1.54\times 10^{-4}$	 & $8.16$&$1.79\times 10^{-2}$&	$2.41\times 10^{-5}$\\
$4000$&	$88.2$&	$4.40\times 10^{-2}$&	$3.18\times 10^{-8}$ &$13.7$&	$3.49\times 10^{-2}$&	$1.22\times 10^{-8}$\\
$8000$&	$107$&	$8.18\times 10^{-2}$&	$1.38\times 10^{-7}$ & $20.5$&	$6.84\times 10^{-2}$&	$2.67\times 10^{-8}$\\
$16000$&$119$&	$1.72\times 10^{-1}$&	$2.29\times 10^{-9}$ & $31.6$&	$1.37\times 10^{-1}$&	$3.35\times 10^{-8}$\\
$32000$&$142$&	$3.19\times 10^{-1}$&	$2.72\times 10^{-9}$ &$51.8$&	$2.73\times 10^{-1}$&	$4.60\times 10^{-8}$\\ 
$64000$&$201$&	$9.47\times 10^{-1}$&	$1.10\times 10^{-8}$ &$94.2$&	$5.69\times 10^{-1}$&	$6.80\times 10^{-8}$\\ \hline \hline

\end{tabular}}
\label{tab:scaling}
\end{table}

 \begin{table}[!b]
\caption{A comparison between the LU decomposition, GMRES, and the fast direct solver when computing $\tau$ and $c$ for the serpentine channel. The timings (in seconds) do not include the initial LU factorization or the fast direct solver's compression of $A^{-1}$.}
\centering
\begin{tabular}{lccc}\hline \hline
$N$ & LU & GMRES & Fast Apply\\ \hline
$500$ & $1.87\times 10^{-2}$ & $23.0$ & $1.11\times 10^{-2}$\\
$1000$ & $5.38\times 10^{-2}$ & $42.3$ & $2.29\times 10^{-2}$\\
$2000$ & $1.89\times 10^{-1}$ & $64.4$ & $4.40\times 10^{-2}$\\
$4000$ & $5.17\times 10^{-1}$ & $91.1$ & $9.06\times 10^{-2}$\\
$8000$ & $-$ & $158$ & $1.64\times 10^{-1}$\\
$16000$ & $-$ & $267$ & $3.70\times 10^{-1}$\\ \hline \hline
\end{tabular}
\label{tab:fastcompare}
\end{table}

\section{Conclusions}
\label{s:concl}
We presented a new algorithm for simulating particulate flows through arbitrary periodic geometries that is comprised mainly of three {\em independent} modules: a periodization scheme for Stokes flow through complex geometries, a fast free-space solver for simulating vesicle flows, and a fast direct solver for the boundary integral equations on the channel. We would like to emphasize that, owing to the linearity of Stokes equations, the periodization scheme can be combined seamlessly with any other existing free-space particulate flow solver implementations (e.g., those for bubbles, drops, rigid-particles, or swimmers). Similarly, other direct solver implementations can be used with equal ease. We showed via numerical experiments that the computational cost of the overall scheme scales linearly with respect to both space and time discretization sizes. 

We are currently working on extending our work on several fronts. To enable higher volume-fraction flow simulations, we are developing a new chunk-based close evaluation scheme to accurately compute the channel-to-vesicle hydrodynamic interactions. In addition, the channel representation will support spatial adaptivity and the channel-vesicle interactions will be treated semi-implicitly. We plan to incorporate islands in the fluid domain by using the standard double-layer formulation as well as extending the scheme to three dimensions in the near future.

\section{Acknowledgements}
GM and SV acknowledge support from NSF under
grants DMS-1224656, DMS-1418964 and DMS-1454010 and a Simons
Collaboration Grant for Mathematicians No. 317933. AB acknowledges
support from NSF under grant DMS-1216656. This research was supported in part through computational resources
and services provided by Advanced Research Computing
at the University of Michigan, Ann Arbor.

\clearpage
\appendix
\section{Arc length correction}
In this section, we present a derivation of the arc length correction formula. The arc length correction is useful for long-time vesicle simulations where the accumulation of errors in the arc length may become significant, leading to elongated vesicles. By placing a small correction term on the right-hand side of the inextensiblity condition, we prevent this accumulation while preserving the membrane's original length with second-order asymptotic convergence. 

We begin by understanding how errors accumulate. Let ${\bf x}^k(\alpha)$, where $\alpha\in[0,2\pi)$, be a parameterization of the membrane $\gamma_k$. The cumulative error on the $k$th time step, denoted by $\mathcal{E}_k$, is given by
\begin{equation}
\mathcal{E}_k=
 \int_{\gamma_k}\! ds-
 \int_{\gamma_0}\! ds.
\label{eq:errordef}
\end{equation}
To find $\mathcal{E}_{k+1}$, we will use the identity
\begin{align}
\|{\bf x}_\alpha^{k+1}\|_2&=\|{\bf x}_\alpha^k+\Delta t {\bf u}_\alpha \|_2\\
&=\left(\left(x_\alpha^k+\Delta t u_\alpha\right)^2+\left(y_\alpha^k+\Delta t v_\alpha\right)^2\right)^\frac{1}{2}\\
&=\left( \left((x_\alpha^k)^2+(y_\alpha^k)^2\right)
+2\Delta t\left( u_\alpha x_\alpha^k+v_\alpha y_\alpha^k \right)
+{\Delta t}^2 \left((u_\alpha)^2+(v_\alpha)^2\right)\right)^\frac{1}{2}\\
&=\left(1+2\Delta t({\bf x}_s^k\cdot {\bf u}_s)
+{\Delta t}^2\|{\bf u}_s\|_2^2\right)^\frac{1}{2}\|{\bf x}_\alpha^k\|_2,
\label{eq:erroran}
\end{align}
where ${\bf x}^k=(x^k,y^k)$ and ${\bf u}=(u,v)$. The error is then
\begin{equation}
\mathcal{E}_{k+1}=\int_{\gamma^k} \left( 1+2\Delta t ({\bf x}_s^k\cdot {\bf u}_s)+{\Delta t}^2 \|{\bf u}_s\|_2^2\right)^\frac{1}{2} \, ds
-\int_{\gamma_0}\! ds.
\label{eq:erroran2}\end{equation}
Observe that setting ${\bf x}_s^k\cdot {\bf u}_s=0$ will cause the arc length to increase monotonically since $\|{\bf u}_s\|_2^2\geq 0$. We now perform a Taylor expansion around $\Delta t=0$ to get
\begin{equation}
\mathcal{E}_{k+1}=\mathcal{E}_k
+\Delta t\int_{\gamma^k}{\bf x}_s^k\cdot {\bf u}_s\, ds
+\frac{1}{2}{\Delta t}^2\int_{\gamma^k}\|{\bf u}_s\|_2^2\, ds
-\frac{1}{2}{\Delta t}^2\int_{\gamma^k}({\bf x}_s^k\cdot {\bf u}_s)^2\, ds
+\mathcal{O}({\Delta t}^3).
\label{eq:arcerror}
\end{equation}
If we set ${\bf x}_s^k\cdot {\bf u}_s=0$, we see that with each time step, the arc length incurs an error on the order of ${\Delta t}^2$. If $n$ is the total number of time steps, we would then expect the cumulative error to scale like $\mathcal{E}_n=\mathcal{O}(n{\Delta t}^2)$. This poses a problem for long-time simulations, where $n$ is very large. To address this, we need to prevent the $\mathcal{E}_k$ term in (\ref{eq:arcerror}) from causing any significant accumulation. A simple remedy is to let
\begin{equation}
{\bf x}_s^k\cdot {\bf u}_s=-\frac{\mathcal{E}_k}{\Delta t \int_{\gamma^k}\! ds}=\frac{\int_{\gamma^0}\! ds-\int_{\gamma^k}\! ds}{\Delta t \int_{\gamma^k}\! ds}.
\label{eq:arccor}
\end{equation}
Each time step, we are still incurring an error on the order of ${\Delta t}^2$. However, we propose that the cumulative error now scales $\mathcal{E}_n=\mathcal{O}({\Delta t}^2)$ and that there exists an upper bound for the relative error that scales ${\Delta t}^2$, independent of $n$.  We only require that $\|{\bf u}_s\|_2$ be bounded. 

To begin, let $L_k=\int_{\gamma^k}\! ds$. Using (\ref{eq:errordef}), (\ref{eq:erroran2}), and (\ref{eq:arccor}), we find that
\begin{equation}
L_{k+1}=\int_{\gamma^k}\left(1+2\frac{L_0-L_k}{L_k}+{\Delta t}^2\|{\bf u}_s\|_2^2\right)^\frac{1}{2}\, ds.
\label{eq:recur}
\end{equation}

\begin{lemma}
Suppose there exists a constant $C$ such that $\|{\bf u}_s\|_2\leq C$ for all time. Then,
\begin{equation}
1-\frac{{\Delta t}^2C^2}{2-{\Delta t}^2C^2}\leq \frac{L_k}{L_0} \leq 1+\frac{{\Delta t}^2 C^2}{2-{\Delta t}^2C^2}
\label{eq:lemmaerror}
\end{equation}
for all $k\geq 0$ when $\Delta t$ is sufficiently small. 
\label{lemma:error}
\end{lemma}\\

\begin{proof}
Clearly, (\ref{eq:lemmaerror}) holds for the base case $L_k=L_0$. Assume the induction hypothesis. We will use the fact that
\begin{equation}
\frac{L_k}{L_0}\left(1+2\frac{L_0-L_k}{L_k}\right)^\frac{1}{2}\leq \frac{L_{k+1}}{L_0} \leq \frac{L_k}{L_0}\left(1+2\frac{L_0-L_k}{L_k}+{\Delta t}^2C^2 \right)^\frac{1}{2}
\end{equation}
to place the desired bounds on $L_{k+1}/L_0$. 

We begin by showing
\begin{equation}
1-\frac{{\Delta t}^2C^2}{2-{\Delta t}^2C^2}\leq\frac{L_k}{L_0}\left(1+2\frac{L_0-L_k}{L_k}\right)^\frac{1}{2}.
\label{eq:longleftcond}
\end{equation}
For simplicity, let $z=L_k/L_0$ and $r={\Delta t}^2C^2$. We need to show that
\begin{equation}
1-\frac{r}{2-r}\leq z\left(\frac{2}{z}-1\right)^\frac{1}{2},\quad \text{for}\quad z\in I_r=\left[1-\frac{r}{2-r},1+\frac{r}{2-r}\right].
\label{eq:leftcond}
\end{equation}
The right-hand side of the inequality in (\ref{eq:leftcond}) is the top half of a circle centered about the point $(1,0)$ with radius 1. The minimum occurs on both endpoints of $I_r$ and is given by
\begin{equation}
\min_{z\in I_r}z\left(\frac{2}{z}-1\right)^\frac{1}{2}=\frac{2}{2-r}\left(1-r\right)^\frac{1}{2}.
\end{equation}
We find that (\ref{eq:leftcond}) holds as long as $r\leq 1$. Therefore, (\ref{eq:longleftcond}) holds as long as $\Delta t\leq 1/C$.

All that remains is to show that
\begin{equation}
\frac{L_k}{L_0}\left(1+2\frac{L_0-L_k}{L_k}+{\Delta t}^2C^2 \right)^\frac{1}{2}\leq 1+\frac{{\Delta t}^2 C^2}{2-{\Delta t}^2C^2},
\label{eq:longrightcond}
\end{equation}
which is equivalent to
\begin{equation}
z\left(\frac{2}{z}-1+r\right)^\frac{1}{2}\leq 1+\frac{r}{2-r},\quad \text{for}\quad z\in I_r.
\label{eq:rightcond}
\end{equation}
The left-hand side of the inequality in (\ref{eq:rightcond}) is the top half of an ellipse centered about $(\frac{1}{1-r},0)$. The maximum occurs on the right endpoint of $I_r$ and is given by
\begin{equation}
\max_{z \in I_r}z\left(\frac{2}{z}-1+r\right)^\frac{1}{2}=1+\frac{r}{2-r}.
\end{equation}
Therefore, (\ref{eq:rightcond}) holds, which completes the proof of Lemma \ref{lemma:error}.
\end{proof}\\

\begin{theorem}
Suppose there exists a constant $C$ such that $\|{\bf u}_s\|_2\leq C$ for all time. Then, the condition
\begin{equation}
{\bf x}_s^k\cdot {\bf u}_s=\frac{L_0-L_k}{\Delta t L_k}, 
\label{eq:ALC} \end{equation}
preserves arc length with relative error
\begin{equation}
\left| \frac{L_{k+1}-L_0}{L_0} \right|\leq \frac{{\Delta t}^2 C^2}{2-{\Delta t}^2 C^2}=\frac{{\Delta t}^2 C^2}{2}+\mathcal{O}({\Delta t}^4)
\end{equation}
for all $k\geq 0$ when $\Delta t\leq 1/C$.
\label{thm:arclength}
\end{theorem}

\newpage
\begin{table}
\caption{Convergence analysis for the arc length and area corrections. To perform the analysis, we evolved a vesicle using time step $\Delta t=0.01/M$ until $t=0.1$. We report the relative errors in the arc length $L$, area $A$, and position ${\bf x}$ with and without the corrections. The subscript $c$ indicates that both the area and arc length corrections were used, and the subscript $nc$ indicates that no corrections were used. The accepted values are labeled with subscript $acc$. We now observe second-order asymptotic convergence for the arc length, which is in agreement with Theorem \ref{thm:arclength}.}
\centering
{\setlength{\extrarowheight}{6pt}
\begin{tabular}{l|ccc|ccc}\hline \hline 
$M$ & $\left|\frac{L_c-L_{acc}}{L_{acc}}\right|$ & $\left|\frac{A_c-A_{acc}}{A_{acc}}\right|$ & $\left\|\frac{{\bf x}_c-{\bf x}_{acc}}{L_{acc}}\right\|_\infty$ & $\left|\frac{L_{nc}-L_{acc}}{L_{acc}}\right|$ & $\left|\frac{A_{nc}-A_{acc}}{A_{acc}}\right|$ & $\left\|\frac{{\bf x}_{nc}-{\bf x}_{acc}}{L_{acc}}\right\|_\infty$\\[6pt] \hline
$1$ & $8.65\times 10^{-5}$ & $1.21\times 10^{-10}$ & $3.02\times 10^{-4}$ & $1.15\times 10^{-3}$ & $3.12\times 10^{-4}$ & $4.78\times 10^{-4}$\\
$2$ & $2.09\times 10^{-5}$ & $7.16\times 10^{-12}$ & $1.55\times 10^{-4}$ & $5.77\times 10^{-4}$ & $1.57\times 10^{-4}$ & $2.38\times 10^{-4}$\\
$4$ & $5.15\times 10^{-6}$ & $4.35\times 10^{-13}$ & $7.80\times 10^{-5}$ & $2.90\times 10^{-4}$ & $7.91\times 10^{-5}$ & $1.17\times 10^{-4}$\\
$8$ & $1.28\times 10^{-6}$ & $2.72\times 10^{-14}$ & $3.90\times 10^{-5}$ & $1.45\times 10^{-4}$ & $3.97\times 10^{-5}$ & $5.68\times 10^{-5}$\\
$16$ & $3.18\times 10^{-7}$ & $3.09\times 10^{-15}$ & $1.93\times 10^{-5}$ & $7.26\times 10^{-5}$ & $1.99\times 10^{-5}$ & $2.65\times 10^{-5}$\\[2.5pt] \hline \hline
\end{tabular}}
\label{tab:correction}
\end{table}

\begin{figure}
\centering
\includegraphics[width=\textwidth]{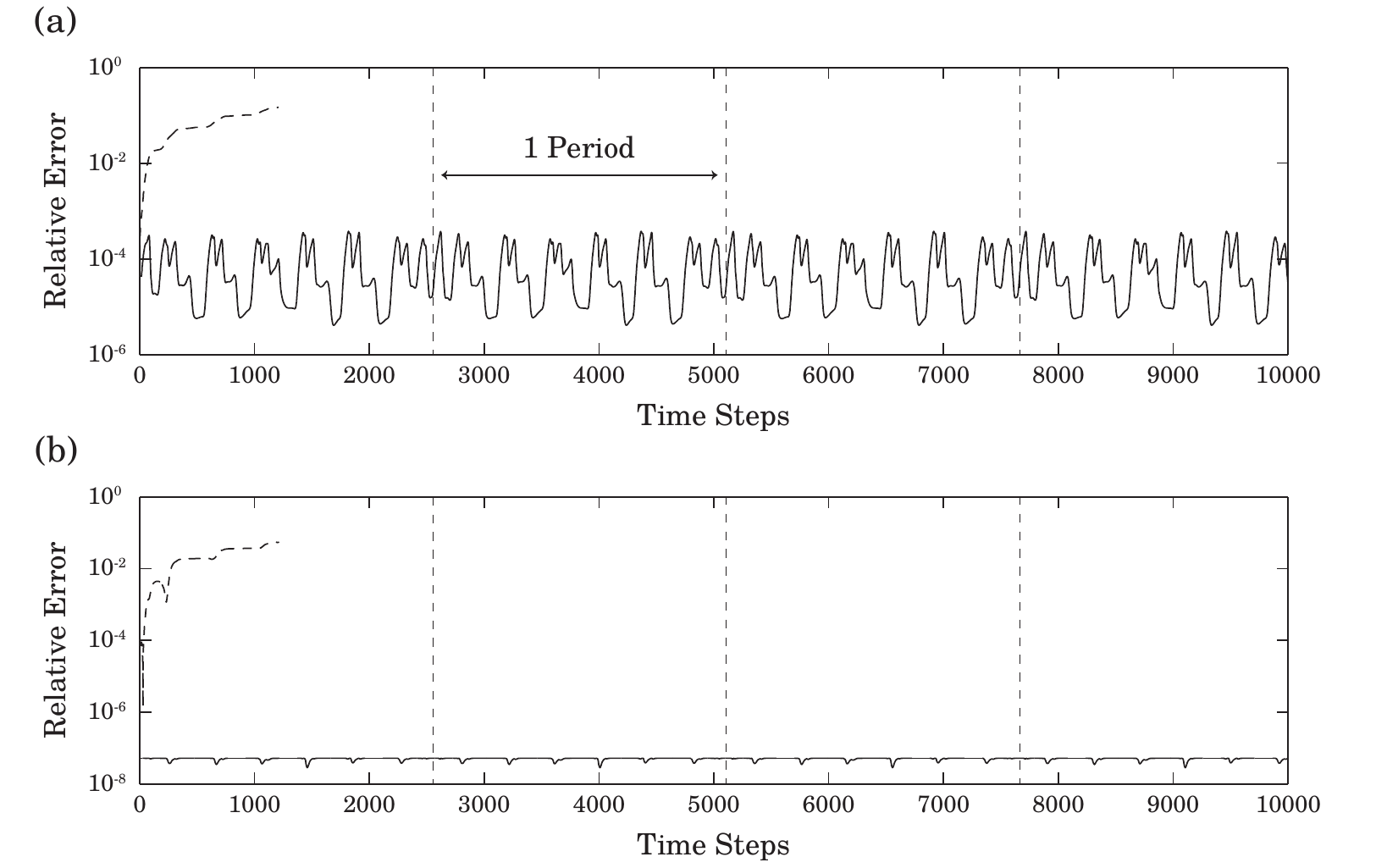}
\caption{Errors without corrections $(--)$ and with both arc length and area corrections $(-)$. The relative error in arc length (a) and area (b) for a single vesicle flowing through the serpentine channel. After around $\sim 1,200$ time steps, the simulation breaks down due to numerical instabilities when no corrections are employed. In the simulation using both corrections, the vesicle successfully passed through the channel several times without incident. In both simulations, the time step was set to $\Delta t=0.01$.}
\end{figure}

\section{Area correction} Although the fluid incompressibility condition is satisfied exactly by the single-layer kernel, errors in the area of vesicles (owing to discretization) accumulate over time and have a compounding effect in the case of long-time simulations. In this section, we discuss a simple and efficient method to correct the area errors whenever the vesicle shapes are updated i.e., at every time step. Let ${\bf x}^k(\alpha)=(x^k(\alpha),y^k(\alpha))$ with $\alpha\in[0,2\pi)$ represent the position of a vesicle's membrane on the $k$th time step. The initial area is given by
\begin{equation}
A_0=\int_0^{2\pi} x^0 y_\alpha^0 \, d\alpha.
\end{equation}
To correct ${\bf x}^k$, we simply add a normal vector ${\bf n}=(y_\alpha^k,-x_\alpha^k)$ scaled by a small unknown constant $c$, which is computed by requiring that the area enclosed by ${\bf x}^k+c{\bf n}$ equals $A_0 $. That is,
\begin{equation}
\int_0^{2\pi} \left( x^k+c y_\alpha^k\right)\left( y^k-cx_\alpha^k \right)_\alpha \, d\alpha =A_0.
\end{equation}
Expanding the integrand gives
\begin{equation}
\int_0^{2\pi} x^ky_\alpha^k \, d\alpha 
-c\int_0^{2\pi}xx_{\alpha \alpha}^k \, d\alpha 
+c\int_0^{2\pi} \left(y_\alpha^k \right)^2 \, d\alpha 
- c^2 \int_0^{2\pi} y_\alpha^k x_{\alpha \alpha}^k \, d\alpha =A_0.
\end{equation}
We take $c$ to be the closest root to zero of this quadratic equation.

\section{Fast reparameterization} The classical approach for reparameterizing evolving geometries in 2D is to introduce an auxiliary tangential velocity in the kinematic condition that helps maintain parameterization quality under some metric (e.g., equispaced in arclength) \cite{hou94}. This approach is very effective and widely used. In our setting, however, in addition to evolving the shape, we need to keep track of the membrane tension from previous time-step. While this can be accomplished by advecting the scalar field (tension) with the auxiliary velocity field, it effects the overall accuracy and stability of our time-stepping scheme. Moreover, our requirements are different: we do not need to maintain exact equi-arclength parameterization but rather a scheme that overcomes the error-amplification due to large shape deformations. For this reason, we take a different approach that avoids an auxiliary advection equation solve at each time step. The scheme only uses Fourier interpolation and redistributes the discretization points on the vesicle membrane so that they are {\em nearly} equispaced. 

Let ${\bf x}(\alpha)$, where $\alpha\in[0,2\pi)$, be a parameterization for a vesicle's membrane $\gamma$. We begin by expressing $ds/d\alpha$ as a Fourier series. That is,
\begin{equation}
\frac{ds}{d\alpha}=\|{\bf x}_\alpha\|_2=\sum_{k=-M/2+1}^{M/2}{C_ke^{ik\alpha}},
\end{equation}
where $M$ is the number of discretization points. Integrating both sides gives
\begin{equation}
s(\alpha)=\int_0^\alpha \|{\bf x}_\beta\|_2 \, d\beta = K+C_0\alpha+\sum_{k\neq 0}{\frac{C_k}{ik}e^{ik\alpha}},
\end{equation}
where
\begin{equation}
K=-\sum_{k\neq 0}{\frac{C_k}{ik}}.
\end{equation}
Our goal is to find a set of discrete points $\{\alpha^*_k\}_{k=0}^{M-1}$ in the parametric domain that corresponds to $\{s_j = s(2\pi j/M), j = 0, \ldots, M-1\}$. We simply use a piecewise linear interpolant of $s(\alpha)$ to determine these points; see Algorithm \ref{algo:computea} for more details. The result is that we obtain points on the curve that are nearly equispaced (up to the accuracy of linear interpolation).  The coordinate positions ${\bf x}$ are then evaluated at $\alpha^*$ via Fourier interpolation. Note that the linear interpolation does not interfere with the overall spectral accuracy of the method (in space) because it is used only to find a new set of points $\alpha^*$. 

\begin{algorithm}
\caption{Compute $\alpha^*$}
\begin{algorithmic}
\State $\alpha^*_0=0$, $j=0$
\For{$k=1:M-1$}
\While{$\frac{s(2\pi)}{M}k\notin[s_j,s_{j+1})$}
\State $j=j+1$
\EndWhile
\State $\alpha^*_k=\frac{2\pi}{M}\left(j+\frac{\frac{s(2\pi)}{M}k-s_j}{s_{j+1}-s_j}\right)$
\EndFor
\end{algorithmic}
\label{algo:computea}
\end{algorithm}

\bibliographystyle{siam}
\bibliography{alex,vesicles,3dVesicles}
\end{document}